\documentclass[a4paper,10t,leqno]{amsart}

\usepackage{amsmath,amssymb,array,euscript,amsfonts}
\usepackage{pifont,graphicx,epsfig}
\usepackage{marginnote}
\usepackage{color}

\def\C{\hbox{\font\dubl=msbm10 scaled 1000 {\dubl C}}}
\def\R{\hbox{\font\dubl=msbm10 scaled 1000 {\dubl R}}}
\def\Q{\hbox{\font\dubl=msbm10 scaled 1000 {\dubl Q}}}
\def\Z{\hbox{\font\dubl=msbm10 scaled 1000 {\dubl Z}}}
\def\N{\hbox{\font\dubl=msbm10 scaled 1000 {\dubl N}}}

\def\sQ{\hbox{\font\dubl=msbm10 scaled 800 {\dubl Q}}}
\def\sR{\hbox{\font\dubl=msbm10 scaled 800 {\dubl R}}}
\def\sZ{\hbox{\font\dubl=msbm10 scaled 800 {\dubl Z}}}
\def\sC{\hbox{\font\dubl=msbm10 scaled 800 {\dubl C}}}

\def\d{\,{\rm{d}}}

\input german.sty

\title[Geometrie der Zahlen]{Geometrie der Zahlen}

\author{Nicola Oswald}

\date{Februar 2016}

\begin{document}

\maketitle

\tableofcontents

\noindent {\small ''G e o m e t r i e \  d e r \  Z a h l e n  habe ich diese Schrift betitelt, weil ich zu den Methoden, die in ihr arithmetische S\"atze liefern, durch r\"aumliche Anschauung gef\"uhrt bin.'' (Hermann Minkowski, 1896)}

\section{Figurierte Zahlen}

Gewisse Verbindungen zwischen Geometrie und Zahlentheorie liegen auf der Hand und treten etliche Male in der ein oder anderen Form in der Geschichte auf. So l\"asst sich einer jeden zusammengesetzten Zahl $n=ab$ mit also nat\"urlichen Zahlen $1<a,b<n$ ein Rechteck mit den Kantenl\"angen $a$ und $b$ zuordnen; im Falle von Primzahlen ist hingegen keine solche Rechteckzerlegung m\"oglich:
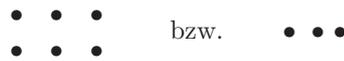
$$
\begin{array}{ccc}
\bullet & \bullet & \bullet \\
\bullet & \bullet & \bullet 
\end{array}\qquad\mbox{bzw.}\qquad \bullet\ \bullet\ \bullet
$$ 
Ein etwas anspruchsvolleres Beispiel figurierter Zahlen bilden die {\it Dreieckszahlen}, welche sich als Partialsumme der Reihe \"uber die nat\"urlichen Zahlen ergeben: 
$$
m=1+2+\ldots+n={\textstyle{1\over 2}}n(n+1).
$$
Der pr\"agnante Ausdruck rechts ergibt sich leicht per Induktion nach $m$ bzw. durch Addition desselben Ausdrucks als Summe mit umgekehrter Summandenfolge:
\begin{eqnarray*}
m &=& 1\ +\qquad 2\quad +\ldots +\ n-1\ +\ n \\
m &=& n\ +\ n-1\quad +\ldots +\quad 2\quad \, +\ \, 1 
\end{eqnarray*}
Auf der rechten Seite bilden wir die Summen untereinanderstehender Zahlen, welche jeweils den Wert $1+n=2+n-1=\ldots n-1+2=n+1$ besitzen; hiervon gibt es $n$ St\"uck, womit im Vergleich der Summe der linken Seiten also $2m=n(n+1)$ bewiesen ist. Diese Beweisidee hatte bereits der Volkssch\"uler Carl Friedrich Gau\ss{} \cite{sartorius}, S. 13. Der Bezug dieser Dreieckszahlen zur Geometrie ergibt sich durch Visualisieren derselben durch \"ubereinandergelagerte Kreise ``$\bullet$'', die insgesamt ein Dreieck bilden; beispielsweise ist $6=1+2+3={1\over 2}\cdot 3\cdot 4$ eine Dreieckszahl:\footnote{Die Leser\_in mag zur \"Ubung einen bildlichen Beweis des folgenden Satzes von Theon von Smyrna aus dem zweiten Jahrhundert vor Beginn unserer Zeitrechnung f\"uhren: {\it  Die Summe zweier aufeinanderfolgender Dreieckszahlen ist eine Quadratzahl} (cf. \cite{deza}).}
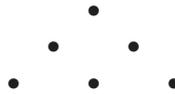
$$
\begin{array}{ccccc}
 & & \bullet & &  \\
 & \bullet &  & \bullet & \\
\bullet & & \bullet & & \bullet
\end{array}
$$
Tats\"achlich hat Gau\ss{} aber etwas viel tiefliegenderes zu Dreieckszahlen bewiesen: In seinem Tagebuch \cite{gaussdiary} findet man als insgesamt 18. Eintrag
$$
\mbox{EYPHKA}!\ \mbox{num}=\Delta+\Delta+\Delta
$$
versehen mit dem Datum 10. Juli 1796. Diese Hieroglyphen sind zu \"ubersetzen als {\it Heureka} (griechisch f\"ur 'ich hab's gefunden', ganz in Anlehnung an Archimedes, der diesen Satz nackt einer Badewanne entspringend bei einer \"ahnlichen Entdeckung gepr\"agt haben soll), {\it jede nat\"urliche Zahl l\"asst sich darstellen als Summe von h\"ochstens drei Dreieckszahlen.} Das Resultat wurde bereits von Pierre de Fermat ohne Beweis ge\"au\ss ert; in Gau\ss ' {\it Disquisitiones} \cite{disqui} findet sich in \S 293 ein erster Beweis.\footnote{Den vielleicht elegantesten Zugang zu diesem wirklich tiefen Satz von Gau\ss{}, liefern die $p$-adischen Zahlen.} Beispielsweise gilt 
$$
2016={\textstyle{1\over 2}}\cdot 63\cdot 64 + 0 + 0= 66 + 465 + 1485= 15 + 231 + 1770 
$$
(um nur einige der zahlreichen derartigen Darstellungen anzugeben). Die weiteren von Fermat ge\"au\ss erten Aussagen \"uber figurierte Zahlen, dass n\"amlich {\it jede nat\"urliche Zahl als Summe von h\"ochstens $n$ so genannten $n$-Eckszahlen dargestellt werden kann}, waren zu Gau\ss ' Tagebuchzeit f\"ur $n\geq 3$ noch unbewiesen, wurden aber 1813 von Augustin Cauchy \cite{cauchy} hergeleitet. Den Spezialfall der Quadratzahlen, wonach jede nat\"urliche Zahl eine Darstellung als Summe von vier Quadraten ganzer Zahlen besitzt, hatte Joseph Louis Lagrange bereits 1770 bewiesen. 
\par

Dreieckszahlen treten \"ubrigens auch bei der so genannten Pellschen Gleichung bzw. dem klassischen Rinderproblem auf. Hierbei handelt es sich um eine von Archimedes an Eratosthenes adressierte Textaufgabe um die Gr\"o\ss e der Herde des Sonnengottes. Die auftretenden acht Teilherden (braune Rinder und dergleichen) unbekannter Gr\"o\ss e  sollen sich als ein Dreieck aufstellen lassen.\footnote{siehe Fowler \cite{fowler}, \S 2.4., f\"ur eine alternative Deutung der zu Eratosthenes' Zeiten wohl unl\"osbaren Aufgabe} 
\par

Die nachstehende Grafik veranschaulicht mit $1+3+\ldots+(2m-1)=m^2$ eine weitere Identit\"at, welche die Summe der ungeraden Zahlen mit den Quadraten in Verbindung setzt:
$$
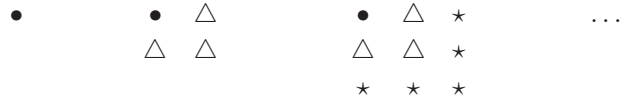
\bullet
\qquad\qquad 
\begin{array}[t]{cc}
\bullet & \triangle \\
\triangle & \triangle 
\end{array}
\qquad\qquad
\begin{array}[t]{ccc}
\bullet & \triangle & \star \\
\triangle & \triangle & \star \\
\star & \star & \star \end{array}
\qquad\qquad \ldots
$$ 
Eine Vielzahl von weiteren Beispielen solcher {\it figurierten Zahlen} finden sich in \cite{deza,rbn}. Tats\"achlich sind figurierte Zahlen jedoch nur am Rande Gegenstand der Geometrie der Zahlen. %Nachstehend ein {\it Beweis ohne Worte} von Deanna B. Haunsperger und Stephen F. Kennedy \cite{hk} und illustriert eine Identit\"at f\"ur Dreieckszahlen $T_k={1\over 2}k(k+1)$, welche zwar seit langer Zeit bekannt ist, aber selten so sch\"on visualisiert wurde. Mehr Beispiele solcher Illustrationen mathematischer Sachverhalte finden sich in Roger B. Nelsens Sammlung \cite{rbn}.
%\begin{figure}[ht]
%\includegraphics[width=12cm]{pww.eps}
%\end{figure}

\section{Das Kreisproblem und seine Verwandten}

\begin{quote}
''Die Grundlage des ganzen Gegenstandes bildet eine eigent\"umliche Untersuchung \"uber die Anzahl aller Combinationen der ganzzahligen Werte, welche zwei unbestimmte ganze Zahlen innerhalb eines vorgeschriebenen Gebietes annehmen. Offenbar kann diese Aufgabe auch unter geometrischer Fassung dargestellt werden, n\"amlich die Anzahl der C o m p l e -\\ x e n \ Z a h l e n zu ermitteln, deren Darstellung innerhalb einer vorgeschreibenen Figur f\"allt.'' 
\end{quote}
schrieb Gau\ss{ } in seinen Untersuchungen {\it \"uber den Zusammenhang zwischen der Anzahl der Klassen, in welche die bin\"aren Formen zweiten Grades zerfallen, und ihrer Determinante} aus dem handschriftlichen Nachlass, niedergeschrieben in der Zeit um 1833\footnote{jedoch basierend auf Ideen um die Jahrhunderwende 1801; zu finden ist diese Schrift im Anhang der {\it Disquisitiones} \cite{disqui}, Ausgabe von 1889, 655-661.}. Der \"Ubergang von zwei (unbestimmten) ganzen Zahlen $a,b$ zu einer komplexen Zahl ist der Darstellung derselben in der komplexen (Gau\ss schen) Zahlenebene geschuldet; heutzutage spricht man bei den $a+ib$ mit $a,b\in\Z$ von den {\it ganzen Gau\ss schen Zahlen}. 
\par

Wir n\"ahern uns dem eigentlichen Themenkreis mit einer einfachen geometrischen Frage: {\it Wie viele Punkte mit ganzzahligen Koordinaten liegen in einem Kreis?} Hier und im Folgenden sprechen wir auch von ganzzahligen Punkten bzw. Gitterpunkten und nennen die Menge $\Z^2$ bestehend aus eben diesen Punkten $(a,b)$ mit ganzzahligen Koordinaten in der euklidischen Ebene ein {\it Gitter}. Sp\"ater werden wir diesen Begriff noch etwas verallgemeinern und pr\"azisieren (im Sinne einer diskreten additiven Untergruppe des $\R^2$ und nicht mit Bezug auf schwedische Gardinen). Zur\"uck zu unserer Frage. Es z\"ahlt 
$$
r(n):=\sharp\{(a,b)\in\Z^2\,:\,a^2+b^2=n\}
$$
f\"ur $n\in\N_0$ einerseits eben die Gitterpunkte $(a,b)\in\Z^2$, welche einen Abstand $\sqrt{n}$ vom Ursprung $(0,0)$ haben, andererseits auch die Anzahl der Darstellungen von $n$ als Summe zweier Quadrate ganzer Zahlen; dies schl\"agt eine Br\"ucke zwischen Geometrie und Zahlentheorie. Nach dem Zweiquadratesatz von Fermat lassen sich alle Primzahlen $p\equiv 1\bmod\,4$ als Summe zweier Quadratzahlen darstellen, so dass also $r(p)>0$ f\"ur ebensolche $p$ gilt; z.B. ist $13=3^2+2^2$. Im Wesentlichen ist jede solche Darstellung (bis auf die aufgrund der Symmetrien m\"oglichen Vorzeichenwechsel und Vertauschung der Summanden) eindeutig, was auf $r(p)=8$ f\"ur diese Primzahlen f\"uhrt. Hingegen ist $r(p)=0$ f\"ur alle Primzahlen $p\equiv 3\bmod\,4$ (weil Quadrate bei Division durch $4$ den Rest $0$ oder $1$ lassen). Dieses recht unterschiedliche Werteverhalten l\"asst sich durch eine Mittelwertbildung besser verstehen: F\"ur $x\geq 0$ ist die Anzahl der ganzzahligen Gitterpunkte in einem Kreis vom Radius $\sqrt{
x}$ um den Ursprung 
$$
R(x):=\sharp\{(a,b)\in\Z^2\,:\,a^2+b^2\leq x\}=\sum_{0\leq n\leq x}r(n).
$$ 
Diese Anzahl von Gitterpunkten l\"asst sich durch die Fl\"ache des Kreises approximieren (verm\"oge eindeutiger Zuordnung eines achsenparallelen Quadrates der Kantenl\"ange eins mit Mittelpunkt in einem jeden solchen Gitterpunkt); also mit Hilfe der Fl\"achenformel $R(x)\sim \pi x$. Mit ein wenig elementarer Geometrie l\"asst sich diese Idee pr\"azisieren:
%\begin{figure}[ht]
%\includegraphics[height=6.85cm]{gauss1803.eps}
%\hspace{0.5cm}
%\includegraphics[height=6.85cm]{geoZah.eps}
%\caption{\footnotesize Der junge Gauss (1803) und seine Idee des Gitterpunktz\"ahlens: Die Anzahl der vom inneren Kreis umschlo\ss enen Gitterpunkte wird approximiert durch die Fl\"achen des inneren und des \"au\ss eren Kreises.}
%\end{figure}
Hierzu fand Gau\ss{}, dass {\it die Anzahl der ganzzahligen Punkte im Inneren eines Kreises vom Radius $\sqrt{x}$ asymptotisch gleich der Kreisfl\"ache ist}, bzw. genauer
$$
\pi \left(\sqrt{x}-{\textstyle{1\over 2}}\sqrt{2}\right)^2< R(x)< \pi\left(\sqrt{x}+{\textstyle{1\over 2}}\sqrt{2}\right)^2
$$
gilt. Diese Ungleichungen lassen sich auch kurz so formulieren\footnote{wobei die Schreibweise $f(x)=O(g(x)$ bedeutet, dass $\vert f(x)\vert \leq Cg(x)$ mit einer absoluten Konstanten $C$ bei $x\to\infty$ gilt} 
$$
R(x)=\pi x+O(x^{1\over 2}),
$$
was unsere Eingangsfrage beantwortet. Insbesondere erweist sich also der Mittelwert von $r(n)$ bei $n\to\infty$ als der Irrationalzahl $\pi$ gleich (womit tats\"achlich auch ein verschwindender Fehlerterm unm\"oglich ist). 
%\begin{figure}[ht]
%\centering
%\includegraphics[height=8cm]{gully.eps}
%\caption{\footnotesize Dieses auf Kreta ans\"assige Muster von Steinchen illustriert die Gau\ss sche Argumentation; Z\"ahlen der 216 Steinchen im Kreis liefert im Vergleich mit dem Quadrat die recht schlechte Approximation ${27\over 8}=3,375$ an die Kreiszahl $\pi=3,141592\ldots$. Tats\"achlich gewann bereits Gau\ss{} durch Abz\"ahlen der Gitterpunkte in immer gr\"o\ss eren Kreisen etliche Dezimalstellen von $\pi$; diese Idee ist verwandt mit so genannten Monte Carlo-Methoden, wie sie bei numerischer Integration verwendet werden.}
%\end{figure} 
Das {\it Kreisproblem} fragt nun nach dem bestm\"oglichen Fehlerterm. Die zur Zeit beste obere Schranke ist $R(x)-\pi x\ll x^{131\over 416}(\log x)^{18\,637\over 8\,320}$ mit ${131\over 416}=0,31490\ldots$ und wurde von Martin Huxley \cite{huxley} erzielt, w\"ahrend andererseits ein Fehlerterm $x^{1\over 4}(\log x)^{1\over 4}(\log\log x)^{1\over 8}$ unm\"oglich ist, wie Kannan Soundararajan \cite{sound} zeigte. Vermutlich ist der Fehlerterm in Wahrheit von der Gr\"o\ss e $x^{{1\over 4}+\epsilon}$.\footnote{Interessant mag in diesem Zusammenhang sein, welche Exponentenjagd sich im Laufe der Forschungst\"atigkeit des letzten Jahrhunderts hierzu stattgefunden hat: So stammen die besten Resultate vor genau einhundert Jahren aus den Federn von Georgi Voronoi \cite{voronoi03} und Godfrey Harold Hardy \cite{hardy}, welche bereits $R(x)-\pi x\ll x^{1\over 3}$ bzw. die Unm\"oglichkeit eines Fehlers $O(x^{1\over 4})$ zeigten, wobei wir hier auf die Angabe von Logarithmuspotenzen der Einfachheit halber verzichten. Hier wurden von bekannten Mathematikern minimale Verbesserungen mit neuen, mitunter recht technischen Absch\"atzungen von Exponentialsummen erzielt; insofern misst der langsame Fortschritt beim Kreisproblem die Schlagfertigkeit der bestehenden Methoden zur Behandlung solcher Aufgaben.} 
\par

Bemerkenswerterweise sind h"oherdimensionale Analoga leichter zu handhaben: Bezeichnet $r_k(n)$ die Anzahl der Darstellungen von $n$ als Summe von $k$ Quadraten, so gilt 
$$
\lim_{x\to\infty}x^{-{d\over 2}}\sum_{n\leq x}r_d(n)={\pi^{d\over 2}\over \Gamma({d\over 2}+1)},
$$
wobei die rechte Seite das Volumen der $d$-dimensionalen Einheitskugel des $\R^d$ ist,\footnote{welche \"ubrigens im Gegensatz zum $d$-dimensionalen Einheitsw\"urfel ein mit $d\to\infty$ gegen null konvergierendes Volumen besitzen(!), wie sich mit Hilfe der Stirlingschen Formel zeigt.} wobei $\Gamma(z)$ die durch $\int_0^\infty t^{z-1}\exp(-t)\d t$ definierte Gamma-Funktion ist. Tats\"achlich ist f\"ur $d\geq 4$ der bestm\"ogliche Fehlerterm bekannt; im Falle $d\geq 5$ ist dieser gegeben durch $O(x^{{d\over 2}-1})$, wie Arnold Walfisz \cite{walfisz} zeigte; die Dimensionen $d=3$ ist hingegen noch nicht vollst\"andig verstanden (siehe auch Fricker \cite{fricker}).\footnote{Warum dieses Gitterpunktproblem im H\"oherdimensionalen einfacher ist, verr\"at die Arithmetik: Die Anzahl $r_k(n)$ der Darstellungen von $n$ als Summe von $k$ Quadraten verh\"alt sich f\"ur $k=2,3$ recht seltsam, mit wachsendem $k$ hingegen wesentlich gleichm\"a\ss iger; siehe hierzu die Anwendungen der nachstehenden Geometrie der Zahlen 
auf Quadratsummen sowie Fricker \cite{frickerhist}.}
\par

Ein verwandtes zweidimensionales Problem besch\"aftigt sich mit Gitterpunkten unterhalb einer Hyperbel. Die Teilerfunktion $n\mapsto d(n):=\sum_{d\mid n}1$ z\"ahlt die positiven Teiler von $n\in\N$. Auch hier oszillieren die Werte von $d(n)$ stark: $d(p)=2$ besteht genau f\"ur prime $p$, w\"ahrend $d(m!)$ mit wachsendem $m$ beliebig gro\ss{} wird. Hier ergibt sich mit der 'Hyperbel-Methode' von Peter Gustav Lejeune Dirichlet \cite{diric42} ein Wachstum \"ahnlich der divergenten harmonischen Reihe, n\"amlich 
$$
\sum_{n\leq x}d(n)=x\log x+(2\gamma-1)x+O(x^{1\over 2}). 
$$
Insofern ist der Mittelwert der Anzahl der verschiedenen Teiler einer nat\"urlichen Zahl $n$ also $\log n$. Das {\it Teilerproblem} fragt nach dem bestm\"oglichen Fehlerterm in dieser Situation. Ersetzt man hier oder beim Kreisproblem Hyperbel oder Kreis durch eine einfach geschlossene {\it glatte} Kurve, so gilt nach einem Resultat von Vojt\v ech Jarn\'\i k\footnote{Jarn\'\i k wurde 1897 in Prag geboren, studierte und promovierte dort, wurde in Prag dar\"uber hinaus auch Professor und verstarb 1970 ebendort. Wichtig f\"ur seine Entwicklung war sicherlich der Aufenthalt beim einflussreichen Edmund Landau in G\"ottingen, der sich ab 1912 selbst eingehend mit dieser Thematik besch\"aftigte. Jarn\'\i k ist bekannt f\"ur seine Arbeiten zu Gitterpunkten in konvexen K\"orpern; in den 1940ern kamen auch Studien zur Geometrie der Zahlen hinzu.} (siehe Steinhaus \cite{jarnik}) \"ubrigens $\vert a-r\vert<\ell$ f\"ur die Anzahl $r$ der von einer beliebigen rektifizierbaren, einfach geschlossenen Kurve der L\"ange $\ell$ umschlossenen Gitterpunkte, deren Inneres Fl\"ache $a$ besitzt.\footnote{Ohne Fehlerterm kommt hingegen folgender erstaunlicher Satz von Georg Alexander Pick \cite{pick} aus: {\it Ist $\Pi\subset \sR^2$ ein konvexes Polygon mit Eckpunkten in $\sZ^2$ und nicht-leerem Inneren, dann gilt $\sharp(\Pi\cap \sZ^2)={\rm{vol}}(\Pi)+{\textstyle{1\over 2}}\sharp(\partial \Pi\cap \sZ^2)+1$, wobei $\partial \Pi$ f\"ur den Rand von $\Pi$ steht.}} Jarn\'\i k lieferte dar\"uber hinaus auch Absch\"atzungen f\"ur die Anzahl von Gitterpunkten auf Kurven \cite{jarnik2}; wichtige Versch\"arfungen f\"ur algebraische Kurven und Anwendungen in der algebraischen arithmetischen Geometrie ergeben sich mit der Arbeit von Bombieri \& Pila \cite{bombieri}.
\smallskip

%\begin{figure}[ht]
%\centering
%\includegraphics[height=5cm]{wald.eps}\qquad\includegraphics[height=5cm]{waldpolya.eps}
%\caption{\footnotesize Manchmal sieht man den Wald vor lauter B\"aumen nicht, mag sich P\'olya gedacht haben. Wenn $r$ gegen\"uber $R$ hinreichend klein ist, l\"asst sich allerdings schon weit blicken. Im Bild rechts sind die schwarz eingef\"arbten 'B\"aume' nicht vollst\"andig sichtbar.}
%\end{figure} 

George P\'olya \cite{polya} stellte sich in seiner nahezu poetisch betitelten Arbeit {\it Zahlentheoretisches und Wahrscheinlichkeitstheoretisches \"uber die Sichtweite im Walde} \cite{polya} die Frage, wenn um jeden Gitterpunkt $(a,b)\neq (0,0)$ in einem Kreis vom Radius $R>0$ um den Ursprung weitere Kreise (B\"aume) vom Radius $r$ geschlagen seien und man selbst im Ursprung $(0,0)$ st\"unde, wie gro\ss{} der kleinste Radius $r$ w\"are, so dass der Blick komplett eingeschr\"ankt w\"are (man also aus dem Wald nicht mehr hinausschauen k\"onnte).\footnote{Wie P\'olya selbst schrieb, geht diese Fragestellung auf seinen Z\"urcher Kollegen und bekannten Gruppentheoretiker Andreas Speiser zur\"uck.} Sind die Radien $r=0$ und $R=\infty$, so ist zwar nicht jeder Gitterpunkt sichtbar, tats\"achlich sind dann nur die $(a,b)\in\Z^2$ sichtbar, f\"ur die $a$ und $b$ teilerfremd sind\footnote{und erste quantitative Absch\"atzungen in dieser Richtung finden sich bereits in der Arbeit \cite{diric42} von einmal mehr Dirichlet. Eine nette Variante hiervon wird von Herzog \& Steward \cite{hs} behandelt: Welche Gitterpunkte sind von einem oder mehreren Gitterpunkten aus sichtbar? Durch Hinzuziehen weiterer Standpunkte erweitert sich nat\"urlich die insgesamt sichtbare Menge der Gitterpunkte; Laison \& Schick \cite{laison} zeigten, dass eine Teilmenge $U\subset \sZ^2$ genau dann unter Zuhilfenahme endlich vieler Standpunkte gesehen werden kann, wenn $U$ f\"ur keine Primzahl $p$ ein perfektes Quadrat modulo $p$ enth\"alt.}, aber sicherlich trifft keine Ursprungsgerade mit irrationaler Steigung einen ganzzahligen Punkt au\ss er den Ursprung. Dieses Problem der F\"orsterei \"uberlassen wir der neugierigen Leser\_in mit dem Hinweis, dass sich P\'olya zur L\"osung der Mathematik des folgenden Paragraphen bedient hat. 

\section{Minkowskis Gitterpunktsatz}

Gau\ss ' Behandlung der Frage nach der Anzahl der Punkte mit ganzzahligen Koordinaten in einem gegebenen Kreis entnehmen wir unmittelbar folgende Verallgemeinerung auf beliebige Gitter: {\it Die Anzahl der Gitterpunkte ${\bf z}\in\Z^2$ in einem Kreis ist asymptotisch gleich dem Volumen des Kreises.} Hermann Minkowski fragte, unter welchen Umst\"anden eine Menge einen Gitterpunkt notwendig enthalten {\it muss.} Dass dies eine f\"ur zahlentheoretische Fragestellungen relevante Beobachtung ist, entnehmen wir bereits dem nachstehenden Zitat 
\begin{quote}
''Wir verdanken Minkowski die fruchtbare Beobachtung, dass gewisse Resultate, welche sich nahezu intutitiv durch Betrachtung von Figuren im $n$-dimensionalen euklidischen Raum ergeben, weitreichende Konsequenzen in verschiedenen Zweigen der Zahlentheorie besitzen.''\footnote{unsere \"Ubersetzung aus dem englischen Original: ''We owe to Minkowski the fertile observation that certain results which can be made almost intuitive by the consideration of figures in $n$-dimensional euclidean space have far-reaching consequences in diverse branchses in number theory.''} 
\end{quote}
von John William Scott Cassels in seinem Klassiker \cite{cassels}, S. 1. Ein Beispiel illustriert diese geometrische Herangehensweise an gewisse zahlentheoretischen Fragestellungen: Zu einer gegebenen reellen Zahl $\alpha$ soll die Existenz einer m\"oglichst guten rationalen Approximation nachgewiesen werden. Das Ma\ss{} der G\"ute einer rationalen N\"aherung ${x\over y}$ ist dabei die Gr\"o\ss e des Nenners. Insofern ist es naheliegend, zu gegebenem ganzzahligen $Q>1$ die beiden Ungleichungen 
\begin{equation}\label{diricl}
\vert y\vert \leq Q\qquad\mbox{und}\qquad \vert y\alpha-x\vert < {1\over Q}
\end{equation}
in ganzen Zahlen $x$ und $y$ erf\"ullen zu wollen; wenn wir vom trivialen Fall $x=y=0$ absehen und also $y\geq 1$ fordern, folgt n\"amlich aus einer solchen L\"osung 
\begin{equation}\label{dirichlet2}
\left\vert \alpha-{x\over y}\right\vert<{1\over yQ}\leq {1\over y^2}.
\end{equation}
Der klassische Approximationssatz von Dirichlet \cite{diric42} besagt, dass zu jedem reellen $\alpha$ und jeder nat\"urlichen Zahl $Q$ es ganze Zahlen $x$ und $y$ gibt, so dass (\ref{dirichlet2}) erf\"ullt ist; wenn dar\"uber hinaus $\alpha$ irrational ist, so gibt es sogar unendlich viele $x$ und $y$ mit dieser Eigenschaft, w\"ahrend es im Falle rationaler $\alpha$ nur endlich viele solcher Paare gibt.\footnote{Dirichlet \cite{diricapprox42} selbst bewies dies und mehrdimensionale Analoga mit einem Schubfachschluss, der ihm oft f\"alschlicherweise als Urheber zugesprochen wird, obwohl diese Schlussweise bereits bei Jean Leurechon im 17 Jahrhundert vorkommt (cf. Heeffer \& Rittaud \cite{heeffer}). Tats\"achlich gibt es neben seinem und dem Minkowskischen Beweis noch weitere Beweise (etwa mittels Kettenbr\"uchen).} F\"ur einen geometrischen Beweis ist also zun\"achst ein Punkt $(x,y)$ mit ganzzahligen Koordinaten verschieden vom Ursprung $(0,0)$ gesucht, der im Inneren des durch die Ungleichungen (\ref{diricl}) beschriebenen achsensymmetrischen Quaders beschrieben ist. Wenn dieser Quader gen\"ugend gro\ss{} ist, dann wird er sicherlich einen solchen Gitterpunkt im Inneren besitzen, und die gew\"unschte Approximation ist gefunden!
\par

Nun gilt es, diesen intuitiv einsichtigen Sachverhalt tats\"achlich zu beweisen. Im Folgenden betrachten wir euklidische Vektorr\"aume $\R^n$ mit stets $n\geq 2$ und dem \"ublichen euklidischen Abstand. Wir denken uns den $\R^n$ gegeben als Menge von Vektoren ${\bf x}=(x_1,\ldots,x_n)$ mit reellen Koordinaten $x_1,\ldots,x_n\in\R$ versehen mit der komponentenweise Addition 
$$
{\bf x}+{\bf y}=(x_1,\ldots,x_n)+(y_1,\ldots,y_n)=(x_1+y_1,\ldots,x_n+y_n)
$$
und der skalaren Multiplikation
$$
\lambda {\bf x}=\lambda (x_1,\ldots,x_n)=(\lambda x_1,\ldots,\lambda x_n)
$$
f\"ur reelle $\lambda$. Der Abstand zwischen zwei Punkten ist dann \"uber den Satz des Pythagoras gegeben durch
$$
d({\bf x},{\bf y})=\sqrt{(x_1-y_1)^2+\ldots+(x_n-y_n)^2}.
$$
Unserem r\"aumlichen Vorstellungsverm\"ogens geschuldet mag man f\"ur die folgenden Ausf\"uhrungen zun\"achst an den zweidimensionalen Fall denken; s\"amtliche Ideen lassen sich unschwer auf den allgemeinen Fall \"ubertragen. 
\par

Minkowskis grundlegende Erkenntnis ist die intuitiv sofort einsichtige wenngleich unpr\"azise Beobachtung, dass eine hinreichend gro\ss e Menge nur dann keinen vom Ursprung ${\bf 0}=(0,\ldots,0)$ verschiedenen Gitterpunkt (mit also ganzzahligen Koordinaten) besitzen kann, wenn sie nicht von spezieller Form ist. Die notwendige Pr\"azisierung erlaubt der Begriff der Konvexit\"at (lat. {\it convexus} f\"ur 'gew\"olbt'). Dieser spielt tats\"achlich eine zentrale Rolle. In der Optik finden sich erste Ideen zur Verwendung gew\"olbter Oberfl\"achen zur Vergr\"o\ss erung bereits in den Schriften des arabischen Universalgelehrten Abu Ali al-Hasan ibn al-Haitham (ca. 965 - 1040); ab dem zw\"olften Jahrhundert fertigten zun\"achst M\"onche aus Beryll erste Linsen (deshalb 'Brille'); die ersten Mikrospkope und Fernrohre wurden gegen Ende des 16. Jahrhunderts von wohl Hans Lippershey und Zacharias Janssen in Holland entwickelt (und kopiert und in eindrucksvoller Weise eingesetzt durch Galileo Galilei). 
\par

Minkowski folgend hei\ss t eine Menge ${\mathcal C}\subset \R^n$ {\it konvex}, falls f\"ur je zwei beliebige ${\bf x}, {\bf y}\in{\mathcal C}$ das verbindende Liniensegment $[{\bf x},{\bf y}]:=\{\ell{\bf x}+(1-\ell){\bf y}\,:\, 0\leq \ell\leq 1\}$ ganz in der Menge ${\mathcal C}$ enthalten ist.\footnote{Bereits bei Archimedes \cite{archi} findet sich eine verwandte Definition konvexer Kurven und Fl\"achen. Tats\"achlich vermeidet Minkowski diesen Begriff zun\"achst und spricht von {\it nirgends konkaven} Mengen; sp\"ater hei\ss t es dann ''Eine nirgends concave Fl\"ache soll als \"uberall convex bezeichnet werden, wenn jede St\"utzebene an die Fl\"ache mit derselben nur einen Punkt gemein hat.'' \cite{minkowski96}, S. 38. In moderner Terminologie werden die Ergebnisse durchweg f\"ur konvexe Mengen formuliert. \"Ubrigens hei\ss t 'konkav' 'einw\"arts gew\"olbt' und stammt ab von dem Lateinischen 'concavus'.} Ferner hei\ss t eine Menge ${\mathcal C}$ {\it symmetrisch}, wenn mit ${\bf x}\in{\mathcal C}$ stets 
auch $-{\bf x}\in{\mathcal C}$ gilt (oder kurz ${\mathcal C}=-{\mathcal C}$).  Nennen wir nun noch eine konvexe Menge ${\mathcal C}$, die nicht in einer Hyperebene enthalten ist, einen {\it konvexer K\"orper}, so k\"onnen wir den ersten Satz von Minkowski formulieren:
\medskip

\noindent {\bf Minkowskis Gitterpunktsatz (1891).} {\it Sei ${\mathcal C}\subset \R^n$ ein symmetrischer konvexer K\"orper mit einem Volumen ${\rm vol}({\mathcal C})>2^n$. Dann enth\"alt ${\mathcal C}$ mindestens einen von ${\bf 0}$ verschiedenen Gitterpunkt ${\bf z}$.}
\medskip

\noindent Eine beschr\"ankte konvexe Menge besitzt stets ein wohldefiniertes endliches Volumen $\mbox{vol}({\mathcal C})$, wie eine approximierende Aussch\"opfung von ${\mathcal C}$ durch Quader zeigt; auch im Falle $n=2$ verwenden wir die Notation ${\rm{vol}}$ f\"ur die Fl\"ache. Besitzt ${\mathcal C}$ unendliches Volumen, so ist die notwendige Ungleichung insbesondere erf\"ullt.
\par

Ein einfacher Beweis hiervon nach Louis J. Mordell \cite{mordell34} basiert auf der Aussch\"opfung des konvexen K\"orpers durch immer kleinere Quader \"ahnlich der Approximation von Integralen durch immer feinere Riemann-Summen. Wir argumentieren im Folgenden nur f\"ur die Ebene (also $n=2$); der allgemeine Fall geht ganz analog. Sei $t$ eine nat\"urliche Zahl. Dann definieren die Gleichungen
$$ 
x_j={2\over t}z_j\qquad \mbox{f\"ur}\quad j=1,2\quad \mbox{mit}\quad z_j\in\Z 
$$
achsenparallele Geraden, welche die Ebene in Quadrate der Fl\"ache $({2\over t})^2$ zerlegen. Bezeichnet nun $N(t)$ die Anzahl der Ecken dieser Quadrate in ${\mathcal C}$, so gilt
$$
\lim_{t\to\infty} \left({2\over t}\right)^2 N(t)={\rm{vol}}({\mathcal C}).
$$
Wegen ${\rm{vol}}({\mathcal C})> 2^2$ ist $N(t)>t^2$ f\"ur hinreichend gro\ss e $t$. 
%\begin{figure}[ht]
%\includegraphics[height=3.3cm]{gauss3.eps}
%\caption{\footnotesize Statt im Kreisproblem den Radius wachsen zu lassen, kann auch der Abstand zwischen den Gitterpunkten verringert werden -- das Verh\"altnis zwischen Gitterpunktanzahl und Fl\"ache der Zellen ist asymptotisch gleich der Kreisfl\"ache. Vielleicht entdeckte Mordell so seinen Beweis des Minkowskischen Gitterpunktsatzes.}
%\end{figure}
Andererseits liefern die P\"archen $(z_1,z_2)$ ganzer Zahlen h\"ochstens $t^2$ verschiedene Paare von Resten bei Division der $z_j$ durch $t$. Demzufolge enth\"alt ${\mathcal C}$ sicherlich zwei verschiedene Punkte
\begin{eqnarray*}
{\sf P}_1:=\left({2z_1^{(1)}\over t},{2z_2^{(1)}\over t}\right),\quad
{\sf P}_2:=\left({2z_1^{(2)}\over t},{2z_2^{(2)}\over t}\right),
\end{eqnarray*}
so dass s\"amtliche $z_j^{(1)}-z_j^{(2)}$ durch $t$ teilbar sind. Aufgrund der Symmetrie von ${\mathcal C}$ ist $-{\sf P}_2\in{\mathcal C}$. Damit besitzt der Mittelpunkt 
$$
{\sf M}=\left({z_1^{(1)}-z_1^{(2)}\over t},{z_2^{(1)}-z_2^{(2)}\over t}\right)
$$
von ${\sf P}_1$ und $-{\sf P}_2$ ganzzahlige Koordinaten, und mit der Konvexit\"at ist ${\sf M}\in{\mathcal C}$. Wegen ${\sf P}_1\neq {\sf P}_2$ ist ${\sf M}\neq {\bf 0}$. Somit ist der Satz (zumindest f\"ur die Ebene) bewiesen.
\medskip

Bislang haben wir den Begriff des Gitters synonym f\"ur die Menge $\Z^n$ verwendet. Allgemein verstehen wir unter einem {\it Gitter} $\Lambda$ in einem euklidischen Raum $V$ eine diskrete additive Untergruppe von $V$, deren Erzeugnis ganz $V$ ist. Im $\R^n$ l\"asst sich mit ein wenig linearer Algebra ein solches Gitter stets darstellen als
$$
\Lambda=\{m_1{\bf z}_1+\ldots +m_n{\bf z}_n\,:\, m_j\in\Z\}=:{\bf z}_1\Z+\ldots+{\bf z}_n\Z
$$
mit gewissen in $\R^n$ linear unabh\"angigen ${\bf z}_1,\ldots,{\bf z}_n$; diese Vektoren ${\bf z}_j$ bilden eine so genannte {\it Basis} des Gitters.\footnote{In diesem Zusammenhang sei erw\"ahnt, dass ganzzahlige Linearkombinationen von \"uber $\sR$ linear abh\"angigen Vektoren keine diskrete Menge bilden.} Der {\it Fundamentalbereich} von $\Lambda$ ist definiert durch
$$
{\mathcal F}=\{\lambda_1{\bf z}_1+\ldots +\lambda_n{\bf z}_n\,:\, 0\leq\lambda_j<1\}\};
$$
die {\it Determinante} des Gitters ist gegeben durch die $n\times n$-Determinante
$$
\det(\Lambda)=\vert\det({\bf z}_1,\ldots,{\bf z}_n)\vert
$$
und entspricht dem Volumen von ${\mathcal F}$ (welches $\neq 0$ ist auf Grund der linearen Unabh\"angigkeit der Spaltenvektoren ${\bf z}_j$). Obwohl die Basis eines Gitters keinesfalls eindeutig bestimmt ist, h\"angt der Wert der Determinante nicht von der Wahl der Basis ab (wie ein wenig lineare Algebra zeigt).
\par

%\begin{figure}[ht]
%\centering
%\includegraphics[height=6.5cm]{cubicspacedivision.eps}\quad \includegraphics[height=6.5cm]{angels.eps}
%\caption{\footnotesize Periodische und aperiodische Parkettierungen sind ein gro\ss es Thema des niederl\"andische K\"unstlers Maurits C. Escher gewesen; oftmals treten in seinen Bildern Gitter explizit auf. Links: Cubic Space Division aus dem Jahr 1952; rechts: Angel-Devil (no. 45) von 1941.}
%\end{figure}

Die bijektive lineare Abbildung
$$
\left(\begin{array}{c} x_1 \\ \vdots \\ x_n \end{array}\right)\mapsto \left(\begin{array}{ccc}
z_{11} & \ldots & z_{1n}\\
\vdots & & \vdots\\
z_{n1} & \ldots & z_{nn}\\
\end{array}\right)\left(\begin{array}{c} x_1 \\ \vdots \\ x_n \end{array}\right)=:\left(\begin{array}{c} y_1 \\ \vdots \\ y_n \end{array}\right)\quad\mbox{mit}\quad {\bf z}_j=\left(\begin{array}{c} z_{1j} \\ \vdots \\ z_{jn}\end{array}\right)
$$
bildet das Standardgitter $\Z^n$ auf $\Lambda$ ab. Dabei werden (symmetrische) konvexe K\"orper des $x$-Raumes in ebensolche im $y$-Raum abgebildet, wobei sich die entsprechenden Volumina proportional um das Verh\"altnis der Volumina der Fundamentalbereiche, also den Faktor $\det(\Lambda)$ ver\"andern.\footnote{Beim strengen mathematischen Beweis geht hier die Transformationsformel der Analysis ein.} W\"are beispielsweise $n=2$ und $\Lambda$ erzeugt durch die Vektoren $({2\atop 0})$ und $({1\atop 1})$, so erg\"abe sich das von diesen beiden Vektoren aufgespannte Parallelogramm als Fundamentalbereich und $\Lambda$ best\"unde aus allen Gitterpunkten von $\Z^2$ mit gerader Koordinatensumme:
$$
\Lambda=\left({2\atop 0}\right)\Z+\left({1\atop 1}\right)\Z=\{\left({a\atop b}\right)\in\Z^2\,:\,a\equiv b\bmod\,2\}.
$$
Ganz allgemein ergibt sich so
\medskip

\noindent {\bf Minkowskis Gitterpunktsatz (allgemeine Version).} {\it Es sei $\Lambda\subset \R^n$ ein Gitter und ${\mathcal C}\subset \R^n$ ein symmetrischer konvexer K\"orper mit einem Volumen ${\rm vol}({\mathcal C})>2^n\det(\Lambda)$. Dann enth\"alt ${\mathcal C}$ mindestens einen vom Ursprung ${\bf 0}$ verschiedenen Gitterpunkt ${\bf z}\in\Lambda$.}
\medskip

\noindent Gilt hingegen ${\rm vol}({\mathcal C})\geq 2^n\det(\Lambda)$, so folgt mit einem Stetigkeitsargument die Existenz eines von ${\bf 0}$ verschiedenen Gitterpunktes in ${\mathcal C}$ oder auf seinem Rand. 
\smallskip

Mit dem Minkowskischen Gitterpunktsatz wird in der Tat das geometrische Konzept der Konvexit\"at f\"ur arithmetische Fragestellungen nutzbar gemacht. Wir illustrieren dies mit folgendem weiteren, auf Harold Davenport \cite{davengazette} zur\"uckgehenden Beispiel aus der klassischen Zahlentheorie, n\"amlich einem Beweis des Zweiquadratesatzes von Fermat\footnote{Fermats Beweis von ca. 1640 ist nicht \"uberliefert, wahrscheinlich argumentierte Fermat mit seiner Abstiegsmethode; siehe hierzu Scharlau \& Opolka \cite{scharl}, S. 8. Den ersten bekannten Beweis erbrachte Leonhard Euler 1749.}: 
%\begin{figure}[ht]
%\centering
%\includegraphics[height=8cm]{13.eps}
%\caption{\footnotesize Ein Beispiel zum Beweis des Zweiquadratesatzes: Das Gitter wird erzeugt durch die Vektoren $(13,0)$ und $(5,1)$. Es ist $q=5$ modulo $13$ eine Quadratwurzel aus $-1$, welche zum Gitterpunkt $(2,3)$ im Kreisinneren f\"uhrt.}
%\end{figure}

Es gilt zu zeigen, dass jede Primzahl $p\equiv 1\bmod\,4$ als Summe zweier Quadrate ganzer Zahlen geschrieben werden kann, wie z.B. $30\,449=100^2+143^2$. F\"ur den Beweis betrachten wir das Gitter 
$$
\Lambda=\left({q\atop 1}\right)\Z+\left({p\atop 0}\right)\Z\subset \R^2\qquad\mbox{mit}\quad q^2\equiv -1\bmod\,p. 
$$
Dass solch ein $q$ tats\"achlich existiert, \"uberlegen wir uns wie folgt: Zu jedem mit $p$ teilerfremden $a$ existiert ein multiplikativ Inverses $a^{-1}$ modulo $p$, welches genau dann von $a$ verschieden ist, wenn $a\not\equiv \pm 1\bmod\,p$ gilt. Entsprechende P\"archenbildung in dem nachstehenden Produkt zeigt 
$$
-1\equiv p!\equiv (-1)^{p-1\over 2}\left(({\textstyle{p-1\over 2}})!\right)^2\bmod\,p,
$$
womit also explizit $q=({p-1\over 2})!$ Gew\"unschtes leistet.\footnote{Dies ist eine Variante des klassischen Satzes von Wilson. Das erste Erg\"anzungsgesetz aus der Theorie der quadratischen Reste liefert unmittelbar die Existenz eines solchen $q$.} Die Fl\"ache des Kreises ${\mathcal C}=\{(x,y)\in\R^2\,:\,x^2+y^2<2p\}$ betr\"agt $2\pi p>4p=2^2\det\Lambda$, womit nach dem Minkowskischen Gitterpunktsatz ein $(a,b)\in{\mathcal C}\cap\Lambda$ verschieden vom Ursprung existiert. Aus der Darstellung
$$
\left({a \atop b}\right)=j\left({q \atop 1}\right)+k\left({p \atop 0}\right)
$$
mit gewissen $k,j\in\Z$ folgt durch Einsetzen
$$
a^2+b^2\equiv j^2(q^2+1)\equiv 0\bmod\,p
$$
nach Wahl von $q$. Also ist $a^2+b^2$ ein Vielfaches von $p$, sicherlich kleiner $2p$, aber nicht $0$. Dies beweist den Zweiquadratesatz.
\smallskip

Ganz \"ahnlich zeigt man den Vierquadratesatz von Lagrange, dass jede nat\"urliche Zahl sich als eine Summe von vier Quadraten darstellen l\"asst. Und tats\"achlich hat Davenport \cite{davengazette} diesen Satz geometrisch bewiesen. Hier wird zun\"achst die Normengleichung f\"ur Quaternionen zum Anlass genommen, dass die Menge der als Summe von vier Quadraten darstellbaren ganzen Zahlen multiplikativ abgeschlossen ist, und man sich beim Nachweis der Aussage des Satzes auf Primzahlen $p\equiv 3\bmod\,4$ zur\"uckziehen kann.\footnote{Der Vierquadratesatz gilt auch in einigen Ganzheitsringen quadratischer Zahlk\"orper; so bewies Fritz G\"otzky \cite{goetzky} dessen G\"ultigkeit in $\sQ(\sqrt{5})$. In der Arbeit von Jesse Ira Deutsch \cite{deutsch} wird ein Beweis mit Minkowskis Geometrie der Zahlen gegeben.} Sogar die noch tiefere Charakterisierung der Zahlen, die sich als Summe von drei Quadraten darstellen lassen, welche zuerst Gau\ss{} fand (und welche seinen Satz \"uber Dreieckszahlen aus dem ersten 
Paragraphen zur Folge hat), gelang Nesmith Ankeny \cite{ankeny} mit dem Minkowskischen Gitterpunktsatz; Vereinfachungen zu seinem Beweis gaben Louis Mordell \cite{mordell58} und Jan W\'ojcik \cite{wojcik}.

\section{Der junge Hermann Minkowski}

Nur selten wurde eine Theorie derma\ss en von einer einzelnen Person in die Welt gesetzt und ihre Fundamente in eine ansprechende Form gebracht, wie es Minkowski mit seiner Geometrie der Zahlen gelang. Tats\"achlich waren die Randbedingungen f\"ur die Sch\"opfung dieser Zahlengeometrie hervorragend. Im Folgenden widmen wir uns den n\"aheren Lebensumst\"anden von Minkowski und stellen die Sichtweise auf Geometrie und Zahlentheorie in seiner Zeit kurz dar.
\smallskip

%\begin{figure}[ht]
%\centering
%\includegraphics[height=4.7cm]{minkokaunas.eps}\quad\includegraphics[height=4.7cm]{stamp.eps}\quad  \includegraphics[height=4.cm,angle=270,origin=c]{grab.eps}
%\caption{\footnotesize Von links nach rechts: Eine Gedenktafel an dem von den Minkowskis besuchten altst\"adtischen Gymnasium; Oskar Minkowski auf einer litauischen Briefmarke verewigt; das Grab der Br\"uder auf dem Berliner Friedhof Heerstra\ss e.}
%\end{figure}

Hermann Minkowski wurde am 22. Juni 1864 im russischen D\"orfchen Aleksotas (heute ein Stadtteil von Kaunas in Litauen) geboren und ging dort kurze Zeit zur Schule, bevor die Familie wegen judenfeindlichen Repressalien\footnote{beispielsweise Studienbeschr\"ankung f\"ur Juden; der \"alteste Bruder Max ging beispielsweise ins benachbarte preussische Insterburg (heute Tschernjachowsk in der russischen Enklave um Kaliningrad) aufs Gymnasium; Bruder Oskar war der erste Jude am Gymnasium in Kaunas). Die begabte j\"ungere Tochter Fanny hingegen besuchte nicht das Gymnasium und ''war Zeit ihres Lebens gegen das Frauenstudium eingestellt'', wie Lily R\"udenberg, die \"alteste Tochter Hermann Minkowskis zu berichten wei\ss{} (\cite{briefe}, S. 13). Frauen hatten zu der Zeit nur an wenigen europ\"aischen Universit\"aten Zugang.} ins nahe preussische K\"onigsberg (heute Kaliningrad in Russland) umzog. Hermann starb unerwartet fr\"uh am 12. Januar 1909 an einer Blinddarmentz\"undung in G\"ottingen. Sein \"alterer Bruder Oskar war ein bedeutender Mediziner, der den Zusammenhang zwischen Blutzucker und der Aktivit\"at der Bauchspeicheldr\"use, was letztlich zur Entwicklung von Insulin f\"uhrte und Oskar den Spitznamen {\it Gro\ss vater des Insulins} einbrachte. Ferner war er Vater des Astrophysikers Rudolph Minkowski und w\"ahrend des ersten Weltkrieges 'Giftgasexperte' auf deutscher Seite. Die Br\"uder Hermann und Oskar liegen in Berlin begraben. 
\par

Bereits als Siebzehnj\"ahriger sorgte Hermann Minkowski mit seiner Arbeit \cite{minko81} zu Darstellungen von nat\"urlichen Zahlen als Summe von f\"unf Quadraten f\"ur Aufsehen und erhielt gemeinsam mit Henry J.S. Smith 1883 den Preis der Pariser Akademie.\footnote{Der ber\"uhmte Physiker Max Born \cite{born}, S. 501, berichtet hierzu: ''So erz\"ahlt Hilbert, da\ss{} Minkowski schon als Student im ersten Semester f\"ur die L\"osung einer mathematischen Aufgabe eine Geldpr\"amie erhielt; er verzichtete aber darauf zugunsten eines armen Mitsch\"ulers und verheimlichte die ganze Sache seiner Familie. Sie ist nur ein Beispiel seiner G\"ute und Bescheidenheit, die alle Menschen erfuhren, welche mit ihm in Ber\"uhrung kamen.''} Kurioserweise war das von der Akademie gestellte Problem bereits 1868 von Smith gel\"ost worden, die L\"osung aber weitgehend unbemerkt geblieben. Minkowski er\"offnete mit seinem Ansatz jedoch das Studium der Formen beliebig vieler Variablen.\footnote{siehe auch Band 1; zu den fr\"uhen Jahren verweisen wir auf Strobl \cite{strobl}.} Zu dieser Zeit nimmt er sein Studium an der hiesigen K\"onigsberger Universit\"at auf. Unter seinen akademischen Lehrern hat wohl Adolf Hurwitz den gr\"o\ss ten Einfluss auf ihn ausge\"ubt. Unter seinen wenigen Kommilitonen befindet sich ein weiterer begnadeter Mathematiker: David Hilbert. Schon bald verbindet die beiden eine enge Freundschaft, die sehr sch\"on in den gesammelten Briefen von Minkowski an Hilbert \cite{briefe} dokumentiert ist.\footnote{Die Briefe von Hilbert an Minkowski sind leider nicht erhalten.} Gemeinsam mit ihrem Lehrer und v\"aterlichen Freund Hurwitz unternehmen die Drei regelm\"a\ss ige Spazierg\"ange; ''Anregungen vermittelten das Mathematische Kolloquium [...], vor allem aber die Spazierg\"ange mit Hurwitz 'nachmittags pr\"azise 5 Uhr nach dem Apfelbaum' '' liest man bei Otto Blumenthal \cite{blumenthal}, dem ersten Doktoranden Hilberts, zu dieser f\"ur das Dreigespann pr\"agenden Zeit. Minkowski studierte zu dieser Zeit auch drei Semester im fernen Berlin und widmete seine Forschungen in dieser Phase seines Schaffens weiterhin den quadratischen Formen; er promovierte 1885 bei Lindemann zu eben diesem Thema. Ferdinand Lindemann (1852-1939) hatte 1882 die Transzendenz von $\pi$ bewiesen und war im darauffolgenden Jahr auf eine Professur in K\"onigsberg berufen worden.\footnote{In Minkowskis Briefen \cite{briefe} an Hilbert finden sich mehrere Stellen, in denen Lindemanns Mathematik deutlich kritisiert wird, so etwa ''Lindemanns Arbeit ist bei n\"aherem Zusehen unter aller Kanone.'' in dem Brief vom 20. September 1901, \cite{briefe}, S. 144.} 
%\begin{figure}[ht]
%\includegraphics[height=5cm]{minkowski.eps}\quad\includegraphics[height=5cm]{albertina.eps}
%\caption{\footnotesize Links: Hermann Minkowski (1864-1909), studierte und promovierte in K\"onigsberg, habilitierte in Bonn; anschlie\ss end hatte Minkowski Professuren in Bonn, K\"onigsberg, Z\"urich und schlie\ss lich G\"ottingen inne. Rechts: Die K\"onigliche Albertus-Universit\"at zu K\"onigsberg.}
%\end{figure}

1887 erlangte Minkowski ein Ruf an die Rheinische Friedrich-Wilhelms-Universit\"at Bonn. In der zugeh\"origen Probevorlesung vom 15. M\"arz 1887 legte er die mit der Rufannahme verbundene Habilitationsschrift {\it R\"aumliche Anschauung und Minima positiv definiter quadratischer Formen}, welche erst mehr als einhundert Jahre sp\"ater 1991 posthum \cite{schwermer} publiziert wurde. In dieser \"au\ss erte Minkowski erste Gedanken zur Geometrie der Zahlen. 
\par

Prominente Vorl\"aufer im Geiste hat Minkowski mit seiner geometrischen Herangehensweise an quadratische Formen in den Gr\"o\ss en Gau\ss{} und Dirichlet. Zun\"achst hatte der Physiker Ludwig August Seeber \cite{seeber} 1831 im Rahmen der Reduktionstheorie quadratischer Formen mit einer monumentalen und undurchsichtigen Arbeit gezeigt, dass  die Determinante einer reduzierten positiv definiten tern\"aren quadratischen Form mindestens halb so gro\ss{} wie das Produkt der Diagonaleintr\"age ist. In seiner umfangreichen Rezension der Seeberschen Arbeit \cite{gauss40} lieferte Gau\ss{} hierf\"ur mit einem geometrischen Ansatz einen ersten vollst\"andigen Beweis; schlie\ss lich gelang Dirichlet 1848 (publiziert jedoch erst 1850 als \cite{dirichlet50}) eine kurze und einfache Darstellung. Er beginnt seine Ausf\"uhrungen gegen\"uber der physikalisch-mathematischen Klasse der K\"oniglich-Preussischen Akademie der Wissenschaften am 31. Juli 1848 in Berlin mit den Worten:\footnote{\cite{dirichlet50}, S. 29}
\begin{quote}
''Indem ich jetzt der Classe das Resultat meiner dahin gerichteten Bem\"uhungen mitzutheilen mir erlaube, glaube ich im Interesse der K\"urze, und wenn ich sagen darf, der Durchsichtigkeit der Darstellung, die geometrische Form beibehalten zu m\"ussen, worin ich die Untersuchung gef\"uhrt habe, der ich die merkw\"urdigen Beziehungen zu Grunde gelegt habe, welche zwischen den quadratischen Formen mit zwei oder drei Elementen und gewissen r\"aumlichen Gebilden Statt finden. '' 
\end{quote}

Hauptergebnis der Minkowskischen Habilitationsschrift ist denn auch die Absch\"atzung des Minimums einer positiv definiten quadratischen Form im Sinne der Vorarbeiten von Gau\ss{} und Dirichlet, diese aber noch weit \"ubertreffend. Joachim Schwermer \cite{schwermer},S. 54, schreibt hierzu: ''Die Frage nach dem Minimum einer positiv definiten quadratischen Form liegt damit an der Quelle des Entstehens der 'Geometrie der Zahlen'. Die Probevorlesung gibt ein fr\"uhes Zeugnis hierf\"ur belegt zugleich, Minkowskis geometrische Ideen zu dieser Zeit und zeigt, da\ss{} er sich der arithmetischen Tragweite seiner \"Uberlegungen schon bewu\ss t war.'' Allerdings berichtet Minkowski erst in einem Brief an Hilbert, datiert 6. November 1889, von ersten Konsequenzen:\footnote{cf. \cite{briefe}, S. 38}
\begin{quote}
''Vielleicht interessirt Sie oder Hurwitz der folgende Satz (den ich auf einer halben Seite beweisen kann): In einer positiven quadratischen Form von der Determinante $D$ mit $n$ ($\geq 2$) Variabeln kann man stets den Variabeln solche ganzzahligen Werthe geben, da\ss{} die Form $<nD^{1\over n}$ ausf\"allt. Hermite hat hier f\"ur den Coefficienten $n$ nur $({4\over 3})^{{1\over 2}(n-1)}$, was offenbar im Allgemeinen eine sehr viel h\"ohere Grenze ist.'' 
\end{quote}
Der angesprochene bedeutende Mathematiker Charles Hermite (1822-1901) hatte in Briefen mit Carl Gustav Jacobi neue Wege in der arithmetischen Analyse quadratischer Formen beschritten. In einem Brief \cite{hermite50} vom 6. August 1845 war ihm die besagte Absch\"atzung mittels eines ausgekl\"ugelten Induktionsbeweises gelungen. Seine Schranke f\"ur dieses erste positive Minimum einer positiv definiten quadratischen Form erlaubt etliche Anwendungen und sie ist f\"ur Minkowski der wesentliche Antrieb, seine Studien zu quadratischen Formen aufzunehmen und letztlich die Geometrie der Zahlen zu entwickeln (siehe hierzu auch Schwermer \cite{schwermer}, S. 51/52).\footnote{Eine Form hei\ss t positiv definit, wenn sie keine negativen Werte annimmt. Indefinite Formen sind wesentlich schwieriger zu behandeln. Beispielsweise korrespondieren die Darstellungen der Null durch die Form $X^2+Y^2-Z^2$ mit den pythagor\"aischen Tripeln.} 
%\begin{figure}[ht]
%\includegraphics[height=4.65cm]{dirichlet.eps}\quad \includegraphics[height=4.65cm]{hermite.eps}\quad \includegraphics[height=4.65cm]{hilbert.eps}
%\caption{\footnotesize Von links nach rechts: Peter Gustav Lejeune Dirichlet (1805-1859), Charles Hermite (1822-1901), David Hilbert (1862-1943)}
%\end{figure}
Angesichts der Minkowskischen Ergebnisse \"au\ss erte Hermite sp\"ater ''Je crois voir la terre promise'' (cf. Hilbert \cite{hilb}, S. 206, bzw. \cite{briefe}, S. 62). Etwas sp\"ater verbesserte Minkowski \cite{minkowski91} sein Resultat noch leicht: 
\medskip

\noindent {\bf Minkowskis Satz \"uber das erste Minimum einer quadratischen Form (1891).} {\it Sei $Q$ eine positiv definite quadratische Form in den Unbekannten $X_1,\ldots,X_n$ mit Determinante $D$, dann existieren ganze Zahlen $x_1,\ldots,x_n$, nicht alle null, so dass}
$$ 
Q(x_1,\ldots,x_n)\leq {\textstyle{4\over \pi}}\Gamma({\textstyle{n\over 2}}+1)^{2\over n}D^{1\over n}.
$$

\noindent Die Gamma-Funktion tritt (wie auch schon beim Kreisproblem) als Ma\ss zahl f\"ur das Volumen der $n$-dimensionalen Kugel auf. F\"ur $n\geq 4$ sticht Minkwoskis Schranke Hermites Absch\"atzung aus. Die {\it Hermite-Konstante} $\gamma_n$ ist definiert als das Infimum\footnote{Leser\_innen, die den Begriff {\it Infimum} nicht kennen, sei mitgeteilt, dass dies die gr\"o\ss te untere Schranke ist; das Minimum einer Menge existiert n\"amlich i.A. nicht. Analog erkl\"art man auch das {\it Supremum} als die kleinste obere Schranke.}, f\"ur die ${\bf x}\in\Z^n\setminus\{{\bf 0}\}$ mit $Q({\bf x})\leq \gamma_n D^{1\over n}$ existiert. Mit der Stirlingschen Formel\footnote{$\Gamma(x+1)=x!\sim \sqrt{2\pi x}(x/e)^x$ bei $x\to\infty$} erweist sich Minkowskis Schranke $\gamma_n\leq(1+o(1)){2n\over \pi e}$ insbesondere f\"ur gro\ss e $n$ als wesentlich kleiner. Sp\"ater erzielte Blichfeldt \cite{blichfeldt29} eine noch bessere Schranke f\"ur die Hermite-Konstante, n\"amlich
$$
\gamma_n\leq {2\over \pi}\Gamma(2+{\textstyle{n\over 2}})^{2\over n}=(1+o(1)){n\over \pi e}.
$$
\"Ubrigens sind nur wenige Werte der Hermite-Konstanten $\gamma_n$ explizit bekannt: $\gamma_2={2\over \sqrt{3}}$ folgt leicht aus der Theorie der bin\"aren quadratischen Formen. Der Wert $\gamma_3=\sqrt[8]{2}$ geb\"uhrt Gau\ss{} \cite{gauss40} in seiner oben zitierten Rezension der Seeberschen Arbeit. Aleksandr Nikolaevich Korkin und Egor Ivanovich Zolotarev \cite{kork,kork2} zeigten $\gamma_4=\sqrt{2}$ sowie $\gamma_5=\sqrt[5]{8}$ noch vor Minkowskis Arbeiten. Mit weiterentwickelten Methoden der Geometrie der Zahlen findet Hans Frederik Blichfeldt in den 1920ern weitere Werte bis einschlie\ss lich $n=8$ (siehe Cassels \cite{cassels}). Dar\"uberhinaus ist noch $\gamma_{24}=4$ bekannt (wof\"ur das Leech-Gitter verantwortlich; siehe \S 11). 
\par

Hilbert schreibt zu Minkowskis Verbesserung der Hermiteschen Schranke f\"ur quadratische Formen in seiner Ged\"achtnisrede auf seinen Freund Minkowski:\footnote{\cite{hilb}, S. 202/3} 
\begin{quote}
''Dieser Beweis eines tiefliegenden zahlentheoretischen Satzes ohne rechnerische Hilfsmittel wesentlich auf Grund einer geometrisch anschaulöichen Betrachtung ist eine Perle Minkowskischer Erfindungskunst. Bei der Verallgemeinerung auf Formen mit $n$ Variablen f\"uhrt der Minkowskische Beweis auf eine nat\"urliche und weit kleinere obere Schranke f\"ur jedes Minimum $M$, als sie bis dahin Hermite gefunden hatte. Noch wichtiger aber als dies war es, da\ss{} der wesentliche Gedanke des Minkowskischen Schlu\ss verfahrens nur die Eigenschaft des Ellipsoides, da\ss{} dasselbe eine konvexe Figur ist und einen Mittelpunkt besitzt, benutzte und daher auf beliebige konvexe Figuren mit Mittelpunkt \"ubertragen werden konnte. Dieser Umstand f\"uhrte Minkowski zum ersten Male zu der Erkenntnis, da\ss{} \"uberhaupt der {\it Begriff des konvexen K\"orpers} ein fundamentaler Begriff in unserer Wissenschaft ist und zu den fruchtbarsten Forschungsmitteln geh\"ort.''
\end{quote}

%\begin{figure}[ht]
%\includegraphics[height=5.4cm]{korkin.eps}\quad \includegraphics[height=5.4cm]{zolotarev.eps}%\quad \includegraphics[height=4cm]{petersburg.eps}
%\caption{\footnotesize Von links nach rechts die Petersburger Mathematiker Aleksandr Nikolaevich Korkin (1837-1908) und Egor Ivanovich Zolotarev (1847-1878). Die beiden bereisten Berlin und insbesondere Paris, wo sie mit Hermite kommunizierten und von ihm mathematisch gepr\"agt wurden. Zolotarev verstarb kurz nach der erfolgreichen gemeinsamen Zusammenarbeit. Deren Arbeit wurde u.a. von dem Mathematiker und Kristallographen Evgraf Stepanovich Fedorov (1853-1919) fortgef\"uhrt. %Ganz rechts das Alexander-Newski-Kloster in St. Petersburg um 1870.}
%\end{figure}

\section{Die Geometrie der Zahlen entwickelt sich}

Die Geometrie fesselte den jungen Minkowski. Bereits 1887 schrieb Minkowski an Hilbert ''Ich bin auch ganz Geometer geworden, und bedaure aus diesem Grunde doppelt, nicht in Ihrem Kreise weilen zu k\"onnen.'' \cite{briefe}, S. 33/34. Das f\"ormliche 'Sie' weicht erst zwei Jahre sp\"ater dem vertrauten 'Du'. Minkowskis Zeit in Bonn hat jedoch auch ihre schlechten Seiten. Hierzu \"au\ss erte Minkowski, wiederum gegen\"uber seinem Freund Hilbert, ''man ist hier ein reiner mathematischer Eskimo'' \"uber seine Zeit in Bonn in einem Brief datiert vom 23. Februar 1893 (siehe \cite{briefe}, S. 51).\footnote{Die Mathematik an der Universit\"at Bonn war seinerzeit dominiert durch den bekannten Analytiker Rudolf Lipschitz, der allerdings zu Minkowskis Zeit gesundheitliche Probleme hatte; Heutzutage erfreut sich Bonn einer Vielzahl hervorragender Mathematik.} Andererseits gelingen Minkowski viele bemerkenswerte Resultate; die Ideen seiner Habilitationsschrift sind sehr fruchtbar.\footnote{Hilbert schreibt in besagter Ged\"achtnisrede auf Minkowski: ''Nunmehr beginnt f\"ur Minkowskis mathematische Produktion die reichste und bedeutendste Epoche; seine bisher auf das spezielle Gebiet der quadratischen Formen gerichteten Untersuchungen erhalten mehr und mehr den gro\ss en Zug ins Allgemeine und gipfeln schlie\ss lich in der Schaffung und dem Ausbau der Lehre, f\"ur die er selbst den treffenden Namen ''{\it Geometrie der Zahlen}''gepr\"agt hat und die er in dem gro\ss artig angelegten Werk gleichen Titels dargestellt hat.'' \cite{hilb}, S. 201}
\par

In dieser Zeit findet Minkowski vielf\"altige Anwendungen seiner geometrischen Methode, etwa in der diophantischen Approximationstheorie und auch in der algebraischen Zahlentheorie \cite{minko90}; den Beweis seines Gitterpunktsatzes (aus \S 3) mit diversen Anwendungen legt Minkowski mit der wichtigen Arbeit \cite{minkowski91} vor.\footnote{Hier findet sich auch die erste Definition von Gittern im $\R^n$; Hlawka w\"urdigte dies mit der Anmerkung: ''Damals hatte der $\R^n$ f\"ur $n>3$ noch lange nicht sein B\"urgerrecht in der Mathematik erworben.'' \cite{hlawka}, S. 10.} Die Vielseitigkeit ist in seiner geometrischen Herangehensweise selbst begr\"undet: Behandelte Minkowski in seinen geometrischen Untersuchungen zun\"achst quadratische Formen, kristallisieren sich nunmehr die Linearformen als zentrales Thema heraus. Eine {\it homogene Linearform} $Y_j$ in den Variablen $X_1,\ldots,X_n$ mit reellen Koeffizienten $\alpha_{ij}$ ist gegeben durch 
$$
Y_j({\bf X})=Y_j(X_1,\ldots,X_n)=\alpha_{j1}X_1+\ldots+\alpha_{jn}X_n.
$$
Unter der Annahme, dass die den Linearformen zugeordnete Determinante $\det(\alpha_{jk})$ nicht verschwindet, bilden die von ganzen Zahlen $x_1,\ldots,x_n$ herr\"uhrenden Tupel $(Y_1(x_1,\ldots,x_n),\ldots,Y_n(x_1,\ldots,x_n))$ ein Gitter $\Lambda\subset\R^n$ mit $\det(\Lambda)=\vert\det(\alpha_{jk})\vert$. Anwenden des Minkowskischen Gitterpunktsatzes auf den durch die Ungleichungen $\vert Y_1({\bf X})\vert\leq\lambda_1\ldots, \vert Y_n({\bf X})\vert\leq\lambda_n$ definierten symmetrisch konvexen K\"orper ${\mathcal C}$ und der Gitterpunktsatz liefert die Existenz ganzer Zahlen $x_1,\ldots,x_n$, nicht alle null, f\"ur die s\"amtliche Linearformen {\it klein} werden. Im Falle verschwindender Determinante ist ${\mathcal C}$ unbeschr\"ankt und eine entsprechende Absch\"atzung gilt trivialerweise:
\medskip

\noindent {\bf Minkowskis Linearformensatz (1891).} {\it Seien Linearformen $Y_1,\ldots,Y_n$ gegeben. F\"ur $\lambda_1\cdot\ldots\cdot\lambda_n\geq \vert\det(\alpha_{jk})\vert$ existiert ${\bf x}\in\Z^n\setminus\{{\bf 0}\}$, so dass }
$$
\vert Y_j({\bf x})\vert\leq \lambda_j\qquad \mbox{f\"ur}\quad j=1,\ldots,n.
$$

\noindent Dies l\"asst sich auf Linearformen mit {\it komplexen} Koeffizienten $\alpha_{jk}$ ausdehnen. Hierf\"ur bildet man Paare konjugiert komplexer Linearformen $Y_k,Y_{k+1}$ und ordnet diesen folgenderma\ss en reelle Paare zu:
$$
{\mathcal Y}_k={\textstyle{1\over \sqrt{2}}}(Y_k+Y_{k+1})\quad\mbox{und}\quad {\mathcal Y}_{k+1}={\textstyle{1\over \sqrt{2}i}}(Y_k-Y_{k+1})
$$
(\"ahnlich den Darstellungen von $\cos$ und $\sin$ durch die Exponentialfunktion). 
\medskip

\noindent {\bf Minkowskis Linearformensatz, komplexe Version (1891).} {\it Gegeben $s$ Paare von Linearformen $Y_1,Y_2,\ldots,Y_{2s-1},Y_{2s}$ mit jeweils komplex konjugierten Koeffizienten $\alpha_{ij}$ und weitere $r=n-2s$ Linearformen $Y_{2s+1},\ldots, Y_n$ mit reellen Koeffizienten. Dann existieren ganze Zahlen $x_1,\ldots,x_n$, nicht alle null, so dass}
$$
\vert Y_j(x_1,\ldots,x_n)\vert\leq \left({\textstyle{2\over \pi}}\right)^{s\over n}\vert\det(\alpha_{jk})\vert^{1\over n}\qquad \mbox{f\"ur}\quad j=1,\ldots,n.
$$

\noindent Mit einem Stetigkeitsargument kann man hier alle bis auf ein Relationszeichen '$\leq$' durch strikte '$<$' austauschen.
\par

Mit den S\"atzen zu Linearformen ergibt sich beispielsweise folgende Versch\"arfung des Dirichletschen Approximationssatzes: {\it Gegeben reelle Zahlen $\alpha_1,\ldots,\alpha_n$, existieren ganze Zahlen $p_1,\ldots,p_n$ und $q\geq 1$ mit}
$$
\left\vert \alpha_j-{p_j\over q}\right\vert<{n\over n+1}{1\over q^{1+{1\over n}}}\qquad\mbox{f\"ur}\quad j=1,\ldots,n.
$$
Dirichlet hatte hier nur den Faktor $1$ statt ${n\over n+1}$. 
\par

''Alle Theoreme hier wiesen {\it einen} Ursprung auf, wir sch\"opften sie aus einer gemeinsamen, sehr durchsichtigen Quelle, die ich als das {\it Prinzip der zentrierten konvexen K\"orper im Zahlengitter} bezeichnen m\"ochte.'' schreibt Minkowski \cite{minkowski07}, S. 234/235, in \"ahnlichem Zusammenhang; er f\"ahrt fort mit den Worten ''Nun sind wir in der Tat eine Strecke Wegs in das Reich der heutigen Zahlentheorie eingedrungen. Wir k\"onnen daran denken, uns auf diesem Boden zu akklimatisieren.'' Und tas\"achlich ergeben sich unmittelbar weitere Anwendungen. Eine solche stellt die algebraische Zahlentheorie bereit: Mit dem komplexen Linearformensatz l\"asst sich die Diskriminante eines Zahlk\"orpers (endliche algebraische Erweiterung von $\Q$) absch\"atzen. Genauer gesagt erzielte Minkowski \cite{minko90} die Ungleichung
$$
\left({4\over \pi}\right)^{r_2}{n!\over n^n}\sqrt{\vert\Delta_K\vert}\geq {\rm{N}}(\mathfrak{a})
$$
f\"ur die Norm eines ganzen Ideals $\mathfrak{a}\neq 0$ im Ganzheistring eines Zahlk\"orpers $K$ vom Grad $n$ \"uber $\Q$ mit $r_2$ Paaren komplexer Einbettungen. Als Konsequenz dieser so genannten Minkowski-Schranke ergibt sich ein einfacher Beweis eines Satzes von Hermite, dass es n\"amlich nur endlich viele Zahlk\"orper mit vorgegebener Diskriminante gibt.\footnote{Auch folgt hieraus die Existenz ganzer Idealen mit kleiner Norm in jeder beliebigen Idealklasse und damit die Endlichkeit der Klassenzahl.} Dar\"uberhinaus zeigt sich, dass es keinen von $\Q$ verschiedenen Zahlk\"orper mit Diskriminante eins gibt \cite{minko96}; konsequenterweise existieren, wie von Kronecker zuvor vermutuet, stets verzweigte Primzahlen \cite{minkowski91}\footnote{ Das Analogon f\"ur Funktionen gilt ebenso: eine algebraische Funktion besitzt stets Verzweigungspunkte.}. Schlie\ss lich gelang Minkowski in diesem Zuge auch noch ein einfacher Beweis des Dirichletschen Einheitensatzes. F\"ur Details verweisen wir auf Neukirch \cite{
neukirch}.

Angesichts der weitreichenden Anwendungen seiner geometrischen Methode kommt Minkowski schlie\ss lich zu dem Schlu\ss{}, seine Theorie durch ein Lehrbuch zu manifestieren. Das Schreiben erweist sich jedoch als schwierig und langwierig: Minkowski ringt um die Form der Darstellung. So entnehmen wir etwa seinem Brief an Hilbert vom 23. Februar 1893:\footnote{\cite{briefe}, S. 50}
\begin{quote}
''Dass ich auch Weihnachten nicht nach K\"onigsberg kam, daran war noch immer mein Buch schuld und die ungem\"utliche Stimmung, in welche ich dadurch versetzt wurde, dass es so langsam fertig wurde. (...) Jetzt steht es mit dem Buche so weit, dass die H\"alfte seit vier Wochen fertig gedruckt ist, und ich das Manuscript der zweiten H\"alfte demn\"achst abgehen lassen will. (...) \"Uberhaupt wollte ich Dich fragen, ob Du vielleicht nicht abgeneigt w\"arest, das Buch nach seinem Erscheinen in den G\"ottinger Anzeigen zu besprechen. Solchen Zwang wie z.B. bei Study w\"urdest Du Dir dabei nicht anzulegen brauchen. Hermite \"ubrigens war von den ihm mitgetheilten Resultaten sehr enthusiasmirt und hat mir bereits den dritten entz\"uckten Brief geschrieben.'' 
\end{quote}
Wir sehen hier Minkowski als selbstbewussten und weitdenkenden Manager seiner Theorie.\footnote{Eine interessante Parallele: Fast zeitgleich zu Minkowskis Niederschrift der Geometrie der Zahlen schrieb Hilbert auf Anfrage der DMV an seinem ber\"uhmten {\it Zahlbericht} zur Darstellung der algebraischen Zahlentheorie inklusive der vielen Erfolge im Zuge der Arbeiten Kroneckers und Dedekinds, zu anfangs \"ubrigens noch gemeinsam mit Minkowski.} Trotzdem erscheint sein Buch \cite{minkowski96} erst 1896; der von Minkowski gew\"ahlte Titel ist nat\"urlich {\it Geometrie der Zahlen} und seine Begr\"undung aus dem Vorwort steht zu Beginn dieses Artikels.
\par

Die Entdeckungen im Zuge seiner Forschungen machen Minkowski zu einem bekannten Zahlentheoretiker; ganz \"ahnlich ergeht es seinem Freund Hilbert. Hilfreich sind die j\"ahrlichen Tagungen der Naturforscher, aber seit der Begr\"undung der Deutschen Mathematiker Vereinigung (DMV) 1890 kommen auch noch deren Jahrestagungen hinzu -- die mathematische Welt ist pl\"otzlich besser vernetzt!\footnote{Dar\"uberhinaus finden zu dieser Zeit die ersten internationalen Mathematiker Kongresse statt, 1893 ein inoffizieller in Chicago, zu dem Minkowski und auch Hilbert k\"urzere Artikel verfassen und im zugeh\"origen Konferenzband \cite{minkowski93} ver\"offentlichen, jedoch ohne selbst nach Amerika zu reisen. 1900 findet in Paris die Weltausstellung und der Internationale Mathematiker Kongress statt, woHilbert seine ber\"uhmte Rede hielt, in der er auf Anraten Minkowskis die wesentlichen Probleme seiner Zeit vorstellte, welche schlie\ss lich (gel\"ost oder ungel\"ost) als die {\it 23 Hilbertschen Probleme} in die Geschichte eingingen (siehe \cite{briefe}, S. 119/120) bzw. \cite{hilberticm}. \"Ahnlich der Geometrie der Zahlen erblickte der Eiffelturm Ende 1887 das Licht der Welt.} 

%\begin{figure}[ht]
%\centering
%\includegraphics[height=4.5cm]{eiffel.eps}\quad \includegraphics[height=4.5cm]{eiffel2.eps}\quad \includegraphics[height=4.5cm]{eiffel3.eps}
%\caption{\footnotesize \"Ahnlich der Geometrie der Zahlen erblickte der Eiffelturm Ende 1887 das Licht der Welt. Seine Fertigstellung zur Weltausstellung im Jahr 1889 zum einhundertj\"ahrigen Geburtstag der franz\"osischen Revolution malte Georges Garen in dem 'Embrasement de la Tour Eiffel pendant l'Exposition universelle de 1889' betitelten Bild. Rechts der w\"ahrend der Weltausstellung 1900 illuminierte Eiffelturm, wie ihn Minkowski und Hilbert beim Besuch des Internationalen Mathematiker Kongresses gesehen haben k\"onnten. Hier hielt \"ubrigens Hilbert seine ber\"uhmte Rede, in der er auf Anraten Minkowskis die wesentlichen Probleme seiner Zeit vorstellte, welche schlie\ss lich (gel\"ost oder ungel\"ost) als die {\it 23 Hilbertschen Probleme} in die Geschichte eingingen (siehe \cite{briefe}, S. 119/120) bzw. \cite{hilberticm}.}
%\end{figure}

\section{Minkowskis Geometrie im Kontext der Zeit}

Geometrie ist eine der mathematischen Disziplinen, welche sich am weitesten in der Geschichte der Menschheit zur\"uck verfolgen l\"asst.\footnote{siehe etwa Scriba \& Schreiber \cite{scriba}} Das Verst\"andnis, was denn unter Geometrie --- buchst\"ablich also der {\it Vermessung der Erde} --- zu verstehen sei, hat sich gerade im 19. Jahrhundert einem gro\ss en Wandel unterzogen. Beginnend mit Thales von Milet im sechsten Jahrhundert vor unserer Zeitrechnung entwickelt sich die Geometrie zun\"achst als die zentrale mathematische Disziplin. Der Bezug zur Zahlentheorie ist indirekt. Zahlen werden meist  als L\"angen, Fl\"achen oder Volumina interpretiert, wie beispielsweie in der Proportionenlehre des Eudoxos mit der geometrischen Wechselwegnahme als Vorl\"aufer des euklidischen Algorithmus. Im dritten Jahrhundert vor unserer Zeitrechnung entwirft Euklid im Rahmen seiner {\it Elemente} eine erste axiomatisch aufgebaute Wissenschaft. Diese euklidische Geometrie wird zweitausend Jahre das Bild von Geometrie pr\"
agen. In einem f\"ur seine Zeit spektakul\"aren Experiment nutzt Eratosthenes etwa zeitgleich geometrische Winkel zur Absch\"atzung des Erdumfangs. Die Aufgaben des Diophant in seiner 'Arithmetica' sind gr\"o\ss tenteils geometrischer Natur.   
\par

Im Anschluss an die klassischen Errungenschaften der griechischen Mathematik ist aus Sicht des europ\"aischen Kulturkreises zun\"achst die analytische Geometrie zu nennen. Diese wurde in der ersten H\"alfte des 17. Jahrhunderts durch Pierre de Fermat und Ren\'e Descartes begr\"undet, Vorabeiten in Form von ersten Koordinaten findet man bereits bei Nioclas Oresme und Fran\c{c}ois Vi\`ete. Bereits bei Fermat werden Tangenten mit Infinitesimalkalk\"ul gebildet; dies wird durch Isaac Newton fortgesetzt und f\"uhrt auf eine Klassifikation der kubischen Kurven. Zeitgleich entwickelt G\'erard Desargues die Anf\"ange der projektiven Geometrie auf der Kunst der Renaissance aufbauend.\footnote{Und in Franken nennt man hier noch den seinerzeit in W\"urzburg und Erlangen ans\"assigen Lehrer Karl Georg Christian von Staudt als weiteren Sch\"opfer der projektiven Geometrie.}
%\begin{figure}[ht]
%\centering
%\includegraphics[height=8.5cm]{duerer.eps}\quad \includegraphics[height=8.5cm]{altend.eps}
%\caption{\footnotesize Der N\"urnberger K\"unstler und Mathematiker Albrecht D\"urer (1471-1528) thematisiert in seinem Buch 'Underweysung der messung mit dem zirckel und richtscheyt' aus dem Jahr 1525 die mathematischen Aspekte der Zentralperspektive (links). Etwa zeitgleich entsteht das erste Landschaftsgem\"alde ohne Menschen, n\"amlich eine Donaulandschaft (rechts) von Albrecht Altdorfer (1480-1538), ein Vorgriff auf die Aufkl\"arung und die nicht mehr zentrale Rolle des Menschen.}
%\end{figure}
Diese Entwicklungen werden Anfang des 19. Jahrhunderts wieder aufgegriffen, zun\"achst durch Gaspard Monge und Jean-Victor Poncelet, sp\"ater August M\"obius und Julius Pl\"ucker; die Erstgenannten untersuchten lineare Transformationen, die Letztgenannten f\"uhrten homogene Koordinaten ein. Daraus entsteht die Invariantentheorie, ein Vorl\"aufer der linearen Algebra, mit Arthur Cayley, James J. Sylvester, Paul Gordan u.a.. Parallel feiert die enumerative Geometrie mit Hermann Hannibal Schubert ihre Hochzeit. Die Entdeckung der nicht-euklidischen Geometrien durch (in alphabetischer Reihenfolge um jeglichen Priorit\"atengejammer aus dem Wege zu gehen) J\'anos Bolyai, Carl Friedrich Gau\ss{} und Nikolai Lobatschewski und letztlich Ferdinand Karl Schweikart liefert ein weiteres gro\ss es Thema der Geometrie des 19. Jahrhundert. 
\par

In Retrospektive florierten zu Lebzeiten Minkowskis insbesondere Analysis und Geometrie. Erste Anwendungen von Methoden dieser Gebiete finden sich in der Zahlentheorie nieder: Beispielsweise in Eulers analytischem Beweis der Unendlichkeit der Menge der Primzahlen oder noch prominenter in Riemanns Programm zur Kl\"arung zentraler Fragen zur Primzahlverteilung von 1859, aber auch die geometrischen Gedanken von Gau\ss{} und Dirichlet im Zusammenhang mit quadratischen Formen (wie oben bereits erl\"autert). 
\par

In diesem Rahmen entwickelt Minkowski seine Ideen. Es ist zu bemerken, dass insbesondere die Entwicklung des passenden Begriffsapparates eine ungemein schwierige Arbeit gewesen sein muss, die ein erhebliches Ma\ss{} an Kreativit\"at erforderte. Minkowski verschmelzt in seiner Geometrie der Zahlen zwei mathematische Begriffe, welche sich zu seiner Zeit gerade verselbstst\"andigen. Da ist zum einen der Begriff des Gitters, welcher sich bereits in den Gau\ss schen Untersuchungen zu quadratischen Formen herauskristallisiert und im Laufe des neunzehnten Jahrhunderts zunehmend an Bedeutung gewinnt. Diesem letztlich algebraischen Begriff steht die geometrische Eigenschaft der Konvexit\"at gegen\"uber. Minkowskis Motivation f\"ur diesen Begriff ergab sich aus dem Problem der Bestimmung des ersten (von null verschiedenen) Minimums einer positiv definiten quadratischen Form. Hier offenbart sich auch die Br\"ucke zwischen Arithmetik und Geometrie, welche Minkowski baut: Die Determinante (bzw. Diskriminante) einer quadratischen Form ist eine geometrische Invariante der hierzu \"aquivalenten Formen. Und beispielsweise im Falle von positiv definiten bin\"aren quadratischen Formen ist die Konvexit\"at der zugeordneten Kegelschnitte die wesentliche Eigenschaft, welche die Existenz eines nicht-trivialen Gitterpunktes bei hinreichender Fl\"ache sicherstellt.    
\par

Karl Hermann Brunn war wohl der Erste, der den Begriff der Konvexit\"at systematisch studierte. Mit seiner Dissertation \cite{brunn87} und seiner Habilitationsschrift \cite{brunn89} aus den Jahren 1887 und 1889 leistete er Grundlegendes zur heutigen Konvexgeometrie (f\"ur einen kurzen Abriss derselben siehe Gruber \cite{gruber90}). Insbesondere die bekannte, oftmals als {\it Ungleichung von Brunn-Minkowski} bezeichnete Absch\"atzung
$$
{\rm{vol}}(\lambda{\mathcal C}_1+(1-\lambda){\mathcal C}_2)^{1\over 3}\geq \lambda {\rm{vol}}({\mathcal C}_1)^{1\over 3}+(1-\lambda){\rm{vol}}({\mathcal C}_2)^{1\over 3}
$$
f\"ur zwei konvexe K\"orper ${\mathcal C}_1, {\mathcal C}_2$ des $\R^3$ und beliebiges reelle $\lambda\in[0,1]$ ist hier zu nennen; dabei ist $\lambda{\mathcal C}_1+(1-\lambda){\mathcal C}_2$ ebenfalls ein konvexer K\"orper. {\it Gilt f\"ur irgendein $\lambda$ Gleichheit, so gilt bereits ${\mathcal C}_2=\kappa{\mathcal C}_1+{\bf x}$ f\"ur ein positives reelles $\kappa$ und ein ${\bf x}\in\R^3$.} Analoga f\"ur h\"ohere Dimensionen bestehen ebenso. Die Brunnsche Beweisf\"uhrung des Gleichheitsfalls wurde von Minkowski, der erst recht sp\"at Brunns Ergebnisse kennen lernte (cf. Kjeldsen \cite{Kjeldsen}, S. 59), kritisiert; beide gaben sp\"ater korrekte Beweise. Minkowski \cite{minko01} fand mit Hilfe der Brunn-Minkowskischen Ungleichung eine neue L\"osung des isoperimetrischen Problems\footnote{dass n\"amlich unter allen geschlossenen Kruven gegebener L\"ange der Kreis die gr\"o\ss te Fl\"ache besitzt. Diese Frage der klassischen griechichschen Mathematik, auch als Didos Problem bekannt, wurde tats\"achlich erst durch die Arbeiten von Jakob Steiner und Friedrich Edler im neunzehnten Jahrhundert gel\"ost. Dank seiner Beliebtheit findet es sich sogar in Leo Tolstois Volkserz\"ahlungen; cf. \cite{hildebrandt}, S. 163.} entgegen einer von Brunn ge\"au\ss erten Erwartung.
%\begin{figure}[ht]
%\centering
%\includegraphics[height=8.5cm]{brunn.eps}\quad \includegraphics[height=8.5cm]{lochkamera.eps}
%\caption{\footnotesize Links: Hermann Karl Brunn (1862-1939). Konvexit\"at tritt insbesondere in der Optik viel fr\"uher auf, etwa bei Lochkameras (camera obscura); rechts eine Skizze zur Lichtbrechung an konvexen und konkaven Linsen von Johannes Zahn im Jahre 1687}
%\end{figure}
Eine detaillierte Analyse, wie sich in Minkowskis Studien (weitestgehend zu quadratischen Formen) aus der geometrischen Intuition \"uber die Konstruktion von so genannten {\it Eichk\"orpern} schlie\ss lich das Konzept der Konvexit\"at herauskristallisiert, liefert Tinne Hoff Kjeldsen \cite{Kjeldsen}.\footnote{In ihrer umfangreichen und sehr lesenswerten Arbeit diskutiert sie auch Minkowskis Arbeiten zu konvexen K\"orpern jenseits der Geometrie der Zahlen.} Sie betont, dass Minkowski in seiner Mathematik dem Trend der Zeit nach Abstraktion und Axiomatisierung gerecht wird\footnote{``This phase also shows Minkowski as a mathematician of the new trend of abstraction and axiomatization that became the hallmark of twentieth century mathematics.'' \cite{Kjeldsen}, S. 86} und ordnet die abstrakten Eichk\"orper als Objekte jenseits unserer reellen Welt ein. In diesem Sinne ist Minkowski in der Tat einer der ersten Vertreter der modernen Mathematik im Sinne der von Jeremy Gray \cite{gray} postulierten mathematischen 
Revolution zur Jahrhundertwende. Sein Freund David Hilbert mit seinen {\it Grundlagen der Geometrie} \cite{hilgeo} von 1899 ist ihm ein Verwandter in diesem Geiste.    
\par

Der weitere wesentliche Begriff in Minkowskis Geometrie neben der Konvexit\"at ist der des Gitters. Tats\"achlich tritt es implizit bereits prominent in Gau\ss{} \cite{gauss40} mit seinen essentiellen Eigenschaften in Erscheinung, allerdings benutzt erst Minkowski die vielseitigen M\"oglichkeiten, die dieser Begriff erlaubt. Strukturmathematisch ist ein Gitter eine diskrete additive Untergruppe eines euklidischen Vektorraums. Diese Definition basiert also auf Konzepten, wie sie gerade zu Minkowskis Zeit ihren Einzug in die Mathematik nahmen. In der Physik, genauer in der Kristallographie beginnend mit Johannes Kepler (auf den wir sp\"ater noch genauer eingehen werden), treten Gitter schon fr\"uhzeitig auf. In der zweiten H\"alfte des 19. Jahrhunderts wurden zunehmend Methoden der jungen Gruppentheorie verwendet. Beispielsweise gelang Arthur Moritz Schoenflies und Jewgraf Stepanowitsch Fjodorow 1890/91 die Beschreibung s\"amtlicher 230 kristallographischen Raumgruppen (cf. Burckhardt \cite{burckhardt}). Auch 
treten Gruppen im Zuge der Dedekindschen Idealtheorie und nicht zuletzt (wenngleich etwas sp\"ater) des Hilbertschen Zahlberichts in dieser Zeit in den Vordergrund zahlentheoretischer Untersuchungen. Insofern ist Minkowski der 
rasanten Entwicklung in verschiedensten Gebieten aufgeschlossen, adaptiert gewisse Sichtweisen und kombiniert sie in genialer Weise zu einem homogenen Ganzen. 
\smallskip

Minkowskis Monographie \cite{minkowski96} kommt sehr gut an. Nach dem kurzen Gastspiel in Bonn kehrt Minkowski 1894 an die K\"onigsberger Universit\"at zur\"uck und beerbt seinen fr\"uheren Lehrer Adolf Hurwitz. Ein Blick auf die Landkarte (im Anhang) zeigt, wie entfernt Minkowski in K\"onigsberg von seinen bisherigen Kontakten war. 1897 erreicht Minkowski ein Ruf an die Polytechnische Hochschule\footnote{die heutige Eidgen\"ossische Hochschule Z\"urich (ETH), an welcher Minkowskis Lehrer Hurwitz seit 1894 t\"atig war; dieser hatte jedoch nicht direkt etwas mit Minkowskis Berufung zu tun (siehe \cite{briefe}, S. 86). Tats\"achlich stand eine Zeit lang auch eine Anstellung an der Universit\"at Z\"urich zur Debatte; Strukturma\ss nahmen standen dem jedoch im Wege, wie man bei Frei \& Stammbach \cite{freist}, S. 46, nachlesen kann.}. 
\par

Auf dem Internationalen Kongress der Mathematiker 1900 in Paris verlautet Hilbert in seiner ber\"uhmten Rede:\footnote{\cite{hilberticm}, S. 259/260}
\begin{quote}
''Die Anwendung der geometrischen Zeichen als strenges Beweismittel setzt die genaue Kenntnis und v\"ollige Beherrschung der Axiome voraus, die jenen Figuren zugrunde liegen, und damit diese geometrischen Figuren dem allgemeinen Schatze mathematischer Zeichen einverleibt werden d\"urfen, ist daher eine strenge axiomatische Untersuchung ihres anschauungsm\"a\ss igen Inhaltes notwendig. Wie man beim Addieren zweier Zahlen nicht unrichtig untereinandersetzen darf, sondern vielmehr erst die Rechnungsregeln, d.h. die Axiome der Arithmetik, das richtige Operieren mit den Ziffern bestimmen, so wird das Operieren mit den geometrischen Zeichen durch die Axiome der geometrischen Begriffe und deren Verkn\"upfung bestimmt.
\par

Die \"Ubereinstimmung zwischen geometrischem und arithmetischem Denken zeigt sich auch darin, da\ss{} wir bei arithmetischen Forschungen ebensowenig wie bei geometrischen Betrachtungen in jedem Augenblicke die Kette der Denkoperationen bis auf die Axiome hin verfolgen; vielmehr wenden wir, zumal bei der ersten Inangriffnahme eines Problems, in der Arithmetik genau wie in der Geometrie zun\"achst ein rasches, unbewu\ss tes, nicht definitiv sicheres Kombinieren an, im Vertrauen auf ein gewisses arithmetisches Gef\"uhl f\"ur die Wirkungsweise der arithmetischen Zeichen, ohne welches wir in der Arithmetik ebensowenig vorw\"arts kommen w\"urden, wie in der Geometrie ohne die geometrischen Einbildungskraft. Als Muster einer mit geometrischen Begriffen und Zeichen in strenger Weise operierenden arithmetischen Theorie nenne ich das Werk von MINKOWSKI ``Geometrie der Zahlen'' (Leipzig 1896).''
\end{quote}
Zuvor hatte Hilbert in seiner Abhandlung \cite{hilbert} sogar gewisse Aspekte von Minkowskis Geometrie der Zahlen als ein Beispiel einer Geometrie nennen, in der das Parallelenaxiom ausgelassen wird.\footnote{siehe hierzu die  k\"urzlich von Volkert neu herausgegebene und kommentierte Hilbertsche Festschrift zu den Grundlagen der Geometrie von 1899 \cite{volkert}, S. 49,76. Umgekehrt war der Enthusiasmus am Hilbertschen Programm der Axiomatisierung der Geometrie nicht zu gro\ss{}. Zassenhaus schrieb ''The arithmetization of geometry which Hilbert tried to establish did not overly excite Minkowski as Minkowski's intuitive geometric ideas stirred up Hilbert's desire to get at the truth without geometric intuition.'' \cite{zasse}, S. 452/453.} 
\smallskip

Schlie\ss lich wird Minkowski 1902 nach G\"ottingen berufen, er ist mittlerweile \"au\ss erst anerkannt. In Retrospektive \"au\ss erte Edmund Hlawka:\footnote{\cite{hlawka}, S. 406}
\begin{quote}
''C. H. Hermite (1822-1901) hat sich \"uber die 1. Lieferung enthusiastisch ge\"au\ss ert. Von R. Fricke (1861-1930) liegt in den Fortschritten der Mathematik eine ausf\"uhrliche Besprechung vor. Ich m\"ochte aber auch die \"Uberlieferung bem\"uhen. Der Physiker Sommerfeld (1868-1951) sagte einmal zu mir, da\ss{} viele Mathematiker nach dem Erscheinen dieses Buches Minkowski \"uber Hilbert (1862-1944) gestellt haben. Er sagte ungef\"ahr, Minkowski habe mehr Ideen gehabt, aber Hilbert w\"are flei\ss iger gewesen. Wenn man die Biographie von Blumenthal (1876-1944) \"uber Hilbert in den gesammelten Werken von Hilbert liest, so findet man ungef\"ahr das gleiche Urteil und das kann nicht ohne Zustimmung von Hilbert erfolgt sein. Heute w\"urde man allerdings dieses Urteil nicht so formulieren.''
\end{quote}

In den Jahresberichten der Deutschen Mathematiker Vereinigung ver\"offentlicht Minkowski anl\"asslich des einhundertsten Geburtstages von Dirichlet eine wahre Lobeshymne auf dessen mathematisches Werk \cite{minkowski05}; dabei offenbart sich einmal mehr, wie geometrische Ideen in Dirichlets Werk Minkowski selbst inspiriert haben: Im Zusammenhang mit einem Fragment in Gau\ss ' Nachlass zur Klassenzahlformel, welche Dirichlet 1837-39 im Zuge seines Beweises f\"ur die Existenz unendlich vieler Primzahlen in einer primen Restklasse entwickelte, erl\"autert Minkowski zun\"achst, dass Gau\ss{} tats\"achlich schon diese Klassenzahlformel kannte, jedoch  Schwierigkeiten beim Nachweis einer Konvergenz hatte, und schreibt ferner:\footnote{in \cite{minkowski05}, S. 456} 
\begin{quote}
''Da\ss{} Dirichlet die fragliche Schwierigkeit \"uberhaupt nicht antrat, liegt daran, da\ss{} er von einem g\"anzlich neuen Ausdruck f\"ur den Fl\"acheninhalt einer Figur ausgehen konnte. Die gew\"ohnliche, ich m\"ochte sagen, {\it mikroskopische} Bestimmung eines Fl\"acheninhalts, welche auch Gau\ss{} auseinandersetzt, besteht darin, auf die zu untersuchende Figur Quadratnetze mit immer engeren und engeren Maschen zu legen und die in die Figur fallenden Maschen zu z\"ahlen. Dirichlet denkt sich ein f\"ur allemal ein bestimmtes, unendliches Quadratnetz fest \"uber die Ebene gebreitet. Die zu untersuchende Figur wird nun von einem festen Anfangspunkt aus kontinuierlich in allen Dimensionen gleichm\"a\ss ig vergr\"o\ss ert, und f\"ur jeden Kreuzungspunkt des Netzes wird angemerkt, bei welchem Vergr\"o\ss erungsverh\"altnis er gerade auf dem Rand der Figur treten w\"urde. Die Summe der reziproken Werte aller so bestimmten Vergr\"o\ss erungszahlen f\"ur alle vorhandenen Netzpunkte w\"urde noch unendlich sein; 
man summiere nun aber bestimmte gleich hohe, etwa $2s^{\mbox{te}}$ Potenzen aller dieser reziproken Werte, so kommt man auf eine durch die urspr\"ungliche Figur v\"ollig charakterisierte Funktion $\zeta(s)$; diese gestattet eine analytische Fortsetzung f\"ur alle komplexen $s$ und weist dann insbesondere f\"ur $s=1$ einen einfachen Pol auf, dessen Cauchysches Residuum genau der Fl\"acheninhalt der Figur wird. Ich m\"ochte vorschlagen, zur leichteren Einb\"urgerung dieser fundamentalen \"Uberlegung in der Analysis die eben beschriebene Formel den Dirichletschen {\it makroskopischen} Ausdruck eines Fl\"acheninhaltes zu nennen.''
\end{quote}
Diese Sprechweise hat sich jedenfalls bis heute nicht durchgesetzt, allerdings greift der Dirichletsche Beweis der analytischen Klassenzahlformel genau diese Idee wunderbar auf. Zudem finden wir Versatzst\"ucke dieses Gedankens auch in Minkowskis Geometrie der Zahlen wieder, welche wir weiter unten genauer er\"ortern wollen.\footnote{Dies bemerkt tats\"achlich bereits der Zeitgenosse Felix Klein \cite{klein}, S. 328, ausdr\"ucklich, wenn er schreibt ''Es findet sich bei ihm eine innere Verwandtschaft mit Dirichletscher Denkweise.'' Und kurz danach gesteht er: ''Ich selbst habe mich seinerzeit darauf beschr\"ankt, gewisse schon bekannte Grundlagen geometrisch klarzustellen, w\"ahrend Minkowski Neues zu finden unternahm. Diese Untersuchungen zeigen deutlich, da\ss{} Geometrie und Zahlentheorie keineswegs einander ausschlie\ss en, sofern man sich in der Geometrie nur entschlie\ss t, diskontinuierliche Objekte zu betrachten.''} 

\section{Sukzessive Minima}

Es lohnt sich, Minkowskis originalen und extrem eleganten Beweis seines Gitterpunktsatzes anzuschauen: Der Einfachheit halber gehen wir vom Gitter $\Lambda=\Z^n$ aus (die Verallgemeinerung auf allgemeine Gitter ist eine sch\"one \"Ubungsaufgabe). Es sei ${\mathcal C}$ eine abgeschlossene Menge ${\mathcal C}\subset \R^n$, die eine Umgebung des Ursprungs enthalte. Selbiges gilt dann auch f\"ur die Menge 
$$
\lambda{\mathcal C}:=\{\lambda{\bf x}\,:\, {\bf x}\in{\mathcal C}\},
$$
wobei $\lambda$ eine positive reelle Zahl sei. Ist $\lambda$ sehr klein, so enth\"alt $\lambda{\mathcal C}$ keinen vom Ursprung verschiedenen Gitterpunkt. In jedem Fall ist die Menge $\lambda{\mathcal C}\cap \Z^n$ endlich. Um nun tats\"achlich einen weiteren Gitterpunkt neben dem Ursprung in einer Menge $\lambda{\mathcal C}$ sicherstellen zu k\"onnen, m\"ussen geometrische Anforderungen an ${\mathcal C}$ selbst gestellt werden. Es ist naheliegend zu fordern, dass keine 'L\"ocher' vorhanden sind (in denen die Gitterpunkte zu liegen kommen k\"onnten). Hier hilft Konvexit\"at. In diesem Fall enth\"alt $\lambda{\mathcal C}$ f\"ur gro\ss e $\lambda$ sicherlich einen von ${\bf 0}$ verschiedenen Gitterpunkt. Wegen der Monotonie
$$
\lambda{\mathcal C}\subset \lambda'{\mathcal C}\qquad\mbox{f\"ur}\quad \lambda<\lambda'
$$
existiert ein kleinstes $\lambda_1$, so dass $\lambda{\mathcal C}$ einen Gitterpunkt ungleich ${\bf 0}$ enth\"alt (und kein $\lambda{\mathcal C}$ f\"ur ein $\lambda<\lambda_1$). Nun betrachten wir die Menge der Translate 
$$
{\bf z}+\lambda{\mathcal C}:=\{{\bf z}+\lambda{\bf x}\,:\,{\bf x}\in{\mathcal C}\}
$$
um Gitterpunkte ${\bf z}\in\Z^n$. Ganz \"ahnlich wie oben zeigt sich, dass f\"ur sehr kleine $\lambda$ diese Translate disjunkt sind, nicht aber f\"ur sehr gro\ss e $\lambda$; insbesondere existiert ein minimales $\lambda_0>0$, so dass
$$
\lambda_0{\mathcal C}\cap ({\bf z}+\lambda_0{\mathcal C})\neq \emptyset\qquad\mbox{f\"ur ein}\quad {\bf z}\neq {\bf 0}. 
$$ 
Wollen wir aus diesem \"Uberschneiden einen Gitterpunkt verschieden vom Ursprung gewinnen, so m\"ussen wir eine weitere geometrische Forderung an unsere Ausgangsmenge stellen: Mit der zus\"atzlichen Vorraussetzung, dass ${\mathcal C}$ symmetrisch ist, folgt f\"ur ein Element ${\bf x}$ beider Mengen, also ${\bf x}\in\lambda_0{\mathcal C}$ und ${\bf x}-{\bf z}\in\lambda_0{\mathcal C}$, dass ebenso ${\bf z}-{\bf x}$ und $-{\bf x}$ in $\lambda_0{\mathcal C}$ enthalten sind. Letzteres zeigt nach Addition von ${\bf z}$, dass zudem ${\bf z}-{\bf x}\in{\bf z}+\lambda_0{\mathcal C}$, womit sowohl ${\bf x}$ als auch ${\bf z}-{\bf x}$ in sowohl $\lambda_0{\mathcal C}$ und ${\bf z}+\lambda_0{\mathcal C}$ liegen. Aufgrund der Konvexit\"at von ${\mathcal C}$ (und selbigem f\"ur $\lambda_0{\mathcal C}$ bzw. dessen Translate) liegt dann auch deren Mittelpunkt
$$
{\textstyle{1\over 2}}({\bf x}\, +\ {\bf z}-{\bf x})={\textstyle{1\over 2}}{\bf z}
$$
in $\lambda_0{\mathcal C}\cap ({\bf z}+\lambda_0{\mathcal C})$. Wir lesen hieraus ${1\over 2}{\bf z}\in \lambda_0{\mathcal C}$ bzw. ${\bf z}\in\lambda_0{\mathcal C}$ und erhalten $2\lambda_0\geq \lambda_1$. Tats\"achlich besteht hier sogar Gleichheit: Enth\"alt n\"amlich $\lambda{\mathcal C}$ einen Gitterpunkt ${\bf z}\neq {\bf 0}$, so ist ${1\over 2}{\bf z}\in{1\over 2}\lambda{\mathcal C}\cap({\bf z}+{1\over 2}\lambda{\mathcal C})$ und speziell f\"ur $\lambda=\lambda_1$ folgt so ${1\over 2}\lambda_1\geq \lambda_0$. 
\par

Weil nun die Translate ${\bf z}+\lambda_0{\mathcal C}$ f\"ur ${\bf z}\in\Z^n$ h\"ochstens Randpunkte gemein haben k\"onnen, muss folglich das Volumen von $\lambda_0{\mathcal C}$ gegen\"uber dem Volumen zwischen den Gitterpunkten beschr\"ankt sein: 
$$
{\rm{vol}}(\lambda_0\mathcal C)\leq 1\qquad\mbox{bzw.}\qquad {\rm{vol}}(\lambda_1\mathcal C)\leq 2^n.
$$
Wenn also kein vom Ursprung verschiedener Gitterpunkt in ${\mathcal C}$ enthalten sein soll, muss $\lambda_1>1$ und damit ${\rm{vol}}(\mathcal C)<2^n$ gelten. Damit ist der Minkowskische Gitterpunktsatz ein zweites Mal bewiesen (und diesmal liefert die leichte Verallgemeinerung auf den Fall allgemeiner Gitter einen Zugang, der weitgehend analytische Hilfsmittel vermeidet). 
\smallskip

Dieser geometrische Zugang Minkowskis suggeriert unmittelbar weitere Fragestellungen. Mit Hilfe der Geometrie der Zahlen lassen sich Gitterpunkte in hinreichen gro\ss en konvexen K\"orpern finden. Dies lieferte beispielsweise eine obere Schranke f\"ur das erste Minimum einer positiv definiten quadratischen Form (s.o.). Es erscheint naheliegend, nach weiteren Gitterpunkten bzw. Minima zu fragen. Entsprechend definiert man zu einem symmetrischen konvexen K\"orper ${\mathcal C}$ und einem gegebenen Gitter $\Lambda\subset \R^n$ die {\it sukzessiven Minima} als
$$
\lambda_j=\inf\{\lambda>0\,:\,\dim(\lambda{\mathcal C}\cap \Lambda)\geq j\}
$$
f\"ur $j=1,\ldots,n$, so ist $\lambda_j$ eine untere Schranke f\"ur alle positiven $\lambda$, so dass $\lambda{\mathcal C}$ mindestens $j$ linear unabh\"angige Gitterpunkte enth\"alt (womit also $\lambda_1$ im Spezialfall der positiv definiten quadratischen Formen mit besagtem Minimum zusammenf\"allt). Die obige Argumentation zeigt bereits
$$
\lambda_1^n\cdot {{\rm{vol}}({\mathcal C})\over \det(\Lambda)}\leq 2^n
$$
Dar\"uberhinaus gelten aber folgende Ungleichungen:
\medskip

\noindent {\bf Minkowskis Satz \"uber sukzessive Minima (1896).}   
$$
{2^n\over n!}\leq \lambda_1\cdot\ldots\cdot \lambda_n\cdot {{\rm{vol}}({\mathcal C})\over \det(\Lambda)}\leq 2^n .
$$

\noindent Der Beweis der oberen Schranke ist recht schwierig (und wir verweisen hier etwa auf Cassels \cite{cassels}). Die Beispiele ${\mathcal C}=\{{\bf x}\in\R^n\,:\,\vert x_1+\ldots+x_n\vert\leq 1\}$ f\"ur die untere Schranke und der W\"urfel ${\mathcal C}=[-1,1]^n$ f\"ur die obere Schranke mit $\Lambda=\Z^n$ zeigen, dass diese Ungleichungen bestm\"oglich sind. Wichtige Verallgemeinerungen f\"ur nicht notwendig symmetrische konvexe K\"orper bewiesen Claude Ambrose Rogers \cite{rogers49} und (unabh\"angig) Claude Chabauty \cite{chabautey49}. Die spektakul\"arste Anwendung in der diophantischen Analysis fand Minkowskis Satz \"uber sukzessive Minima wohl in der scharfen Form des Siegelschen Lemmas, welche Enrico Bombieri \& Jeffrey Vaaler \cite{vaaler} bewiesen.   
\par

Die sukzessiven Minima sind jedoch kein vollst\"andig gekl\"artes Terrain. Zu einem symmetrischen konvexen K\"orper ${\mathcal C}$ wird die zugeh\"orige {\it kritische Determinante} $\Delta({\mathcal C})$ als das Infimum aller Gitterdeterminanten $\det(\Lambda)$ zu Gittern $\Lambda$ erkl\"art, welche keinen von ${\bf 0}$ verschiedenen Gitterpunkt in ${\mathcal C}$ besitzen. Minkowski gelang der Nachweis der Ungleichung
$$
\lambda_1\cdot\ldots\cdot \lambda_n\cdot \Delta({\mathcal C})\leq \det(\Lambda)
$$
f\"ur den Fall $n=2$ und $n$-dimensionaler Ellipsoide; Alan C. Woods \cite{woods} zeigte selbiges f\"ur $n=3$, der allgemeine Fall ist aber bis heute offen. 
\par

Ein verwandter weiterer Aspekt, den Minkowskis originaler Beweis des Gitterpunktsatzes anregt, betrifft die Frage, wie dicht sich gegebene konvexe K\"orper im Raum packen lassen. Beispielsweise kann man mit den Translationen ${\mathcal F}+{\bf z}$ des Fundamentalbereichs des Gitters um Gitterpunkte den gesamten Raum l\"uckenlos und ohne \"Uberdeckung pflastern. Hieran ankn\"upfende Untersuchungen werden wir im Rahmen unserer Analyse der Rezeption der Minkowskischen Arbeiten vorstellen.  

\section{Minkowskis Raum-Zeit}

In seinem Vortrag anl\"asslich des Internationalen Mathematiker Kongresses 1904 in Heidelberg stellt Minkowski \cite{minkowski04} sein Gebiet der {\it Geometrie der Zahlen} als ''introduction des variables continues dans la th\'eorie des nombres'' mit Worten von Charles Hermite vor; er pr\"azisiert dies mit ''Einige hervorstechende Probleme darin betreffen die Absch\"atzung der kleinsten Betr\"age kontinuierlich ver\"anderlicher Ausdr\"ucke f\"ur ganzzahlige Werte der Variablen.'' \cite{minkowski04}, S. 164. Seinen Gitterpunktsatz bezeichnet er dabei selbstbewusst als ''das Fundamentaltheorem der Geometrie der Zahlen (...) weil es fast in jede Untersuchung auf diesem Gebiete hineinspielt.'' \cite{minkowski04}, S. 164. In seinem Vortrag und dem zugeh\"origen Artikel stellt Minkowski das gesamte Spektrum von Anwendungen seiner Geometrie der Zahlen in eindrucksvoller Weise dar. So stellt er heraus, dass dieser geometrische Ansatz ''zu einfachen Beweisen der Dirichletschen S\"atze \"uber die Einheiten in den 
algebraischen Zahlk\"orpern''\footnote{\cite{minkowski04}, S. 168} f\"uhrt, und tats\"achlich findet sich in der modernen Literatur zur algebraischen Zahlentheorie nahezu ausnahmslos der Minkowskische Beweis des Dirichletschen Einheitensatzes\footnote{der da besagt, dass die Einheitengruppe des Ganzheitsringes eines Zahlk\"orpers mit $R_1$ reellen und $r_2$ Paaren komplex konjugierter Einbettungen in $\sC$ eine abelsche Gruppe vom Rang $r_1+r_2-1$ ist. Im Falle reell-quadratischer Erweiterungen $\sQ(\sqrt{d})$ von $\sQ$ liefert die Minimall\"osung der Pellschen Gleichung $X^2-dY^2=1$ ein erzeugendes Element.}. Ferner spricht Minkowski auch mehrdimensionale Kettenbr\"uche an\footnote{wof\"ur wir auf Erd\H{o}s et al. \cite{erdoes} verweisen; hierzu sei angemerkt, dass etwa zeitgleich Klein \cite{kleinkb} geometrisch verwandte Studien durchf\"uhrte.} und liefert hierf\"ur eine geometrische Herangehensweise. Und seine Fragezeichen-Funktion $?(x)$ mit ihren seltsamen analytischen Eigenschaften wird erw\"ahnt; sie 
ist bis heute Gegenstand tiefsinniger Untersuchungen im Kontext der Arithmetik von Kettenbr\"uchen. Res\"umierend kann man sagen, dass Minkowskis Geometrie der Zahlen in k\"urzester Zeit ein weitreichendes Spektrum an Anwendungen gefunden hat.

%\begin{figure}[ht]
%\centering
%\includegraphics[height=8cm]{minkobild.eps}
%\caption{\footnotesize Dem anl\"asslich Minkowskis Vortrag auf dem Internationalen Mathematiker Kongress 1904 in Heidelberg ver\"offentlichen Artikel \cite{minkowski04} beigef\"ugtes Faltblatt zur Illustration.}
%\end{figure}

Drei Jahre sp\"ater verfasst Minkowski mit seinem Lehrbuch {\it Diophantische Approximationen} \cite{minkowski07} eine leichter lesbare Einf\"uhrung in das von ihm durch seine Geometrie der Zahlen bereicherte Gebiet; bis zu Jurjen Koksmas Lehrbuch gleichen Namens \cite{koksma} aus dem Jahr 1936 f\"allt ihm die Rolle der Standardreferenz zu diophantischen Fragestellungen zu. Innovativ in Minkowskis Lehrbuch ist auch die zentrale Rolle, welche er der Arithmetik von Zahlk\"orpern zuweist, ein Thema, welches zu dieser Zeit besonders aktuell war (nicht zuletzt dank des {\it Zahlbericht}s seines Freundes Hilbert). 
\smallskip

Grob betrachtet mag man Minkowskis Schaffen in drei Phasen einteilen: die zahlentheoretische Fr\"uhphase zu quadratischen Formen, die Schaffenszeit der 1890er beginnend mit der Berufung nach Bonn und dem Verfassen der Geometrie der Zahlen, und schlie\ss lich seine Physikphase ab 1900. Tats\"achlich ist die Physik dieser Epoche durchdrungen von Geometrie.\footnote{In diesem Kontext ein weiteres Zitat von Born: ''Der gr\"o\ss te Teil von Minkowskis Arbeiten liegt auf dem Felde der Zahlentheorie, und davon zu berichten, bin ich nicht befugt. Ich habe nur in einem seiner arithmetischen Werke gelesen, das den Titel 'Diophantische Approximationen' tr\"agt und aus seinen Vorlesungen im Winter 1903/04 hervorgegangen war, gerade ehe ich nach G\"ottingen kam. Darin werden wie in Minkowskis gr\"o\ss erem Werk 'Geometrie der Zahlen' arithmetische S\"atze aus geometrischen Betrachtungen gewonnen, und zwar hier insbesondere mit Hilfe des Netzes von Punkten mit ganzzahligen Koordinaten in der Ebene. Diese Einf\"uhrung in 
das Zahlengitter war mir sp\"ater bei der dynamischen Theorie der Kristallgitter von gewissem Nutzen.'' \cite{born}, S. 41. Tats\"achlich ver\"offentlich Born \cite{born15} (gewidmet seinem Freunde David Hilbert) 1915 seine vielbeachtete Monographie {\it Dynamik der Kristallgitter}, welche zusammen mit seinem Artikel zur Atomtheorie des festen Zustandes von 1923, die Gitterdynamik in einheitlicher Form darstellt, was heutzutage als ein wichtiger Schritt zur Begr\"undung der Festk\"orperphysik angesehen wird. Den Physik-Nobelpreis erhielt Born 1954 jedoch f\"ur seine Arbeiten zur statistischen Deutung der Quantenmechanik.} Sein zunehmendes Interesse an physikalischen Fragestellungen fand 1907 einen H\"ohepunkt, als Minkowski realisierte, wie sich die Theorien von Lorentz und Einstein mathematisch in einem nicht-euklidischen Raum erkl\"aren lassen. Nach Einstein sind Messwerte f\"ur Zeit und Raum relativ zum Beobachtenden, Minkowski erkannte, dass eine bestimmte Kombination derselben unabh\"angig vom 
Beobachtenden ist. Die Konsequenz ist sein Raum-Zeit-Kontinuum, welches sp\"ater Einstein bei der Entwicklung dessen allgemeiner Relativit\"atstheorie hilfreich war. Auch hier finden sich quadratische Formen (in Gestalt der so genannten Minkowski-Metrik) und geometrische Gedanken (die Verbindung von Raum und Zeit zu einer Raumzeit). Max Born w\"urdigt die Beitr\"age der Beteiligten wie folgt:\footnote{\cite{born}, S. 504}
\begin{quote}
''Im ganzen kann man sagen, da\ss{} die spezielle Relativit\"atstheorie nicht das Werk eines Mannes, sondern durch Zusammenwirken einer Gruppe gro\ss er Forscher, LORENTZ, POINCARE, EINSTEIN, MINKOWSKI, entstanden ist. Da\ss{} gew\"ohnlich EINSTEINs Name allein genannt wird, hat gleichwohl eine gewisse Berechtigung, weil ja die spezielle Theorie nur der erste Schritt war zu der allgemeinen, welche die Gravitation mit umfa\ss te und dadurch das gesamte Werk NEWTONs revolutionierte. Die allgemeine Relativit\"atstheorie aber ist EINSTEINs ausschlie\ss liche Leistung. Sie beruht auf der Verbindung der Minkowskischen Weltgeometrie und den tiefen Gedanken \"uber gekr\"ummte R\"aume, die lange vorher von BERNHARD RIEMANN entwickelt und von CHRISTOFFEL, RICCI und LEVI-CIVIT\'A weiter gef\"uhrt worden waren. Auch die allgemeine Relativit\"atstheorie ist somit ohne Minkowskis Arbeit undenkbar, und darum ist es nicht ohne Reiz zu fragen, was EINSTEIN von MINKOWSKI gehalten hat. In der ersten Zeit, um 1909, da ich 
EINSTEIN kennenlernte, war er ziemlich ablehnend und sah in MINKOWSKIs Arbeit nicht viel mehr als \"uberfl\"ussiges mathematisches Beiwerk. Aber das \"anderte sich schnell, als er tiefer in das Problem der allgemeinen Relativit\"at eindrang, wo gerade MINKOWSKIs mathemathische Methoden wesentlich wurden.'' 
\end{quote}
Interessant ist, dass Einstein Minkowski aus dessen Z\"uricher Zeit sehr wohl bekannt war; hierzu schreibt Born \"uber Minkowski:\footnote{\cite{born}, S. 502} 
\begin{quote}
''Unter seinen Sch\"ulern war einer, dessen Name kurze Zeit sp\"ater mit dem seinen viel genannt werden sollte, als die spezielle Relativit\"attstheorie die Gem\"uter bewegte, ALBERT EINSTEIN. Aber er ist Minkowski keineswegs besonders aufgefallen. Als ich sp\"ater, 1909, Minkowskis Mitarbeiter an Problemen der Relativit\"atstheorie geworden war, hat er mir einmal gesagt: ''Ach, der Einstein, der schw\"anzte immer die Vorlesungen -- dem h\"atte ich das gar nicht zugetraut.'''' 
\end{quote}

%\begin{figure}[ht]
%\centering
%\includegraphics[height=5cm]{relativity.eps}\quad \includegraphics[height=5cm]{picasso.eps}
%\caption{\footnotesize Links ein sich mit 90 Prozent der Lichtgeschwindigkeit n\"aherndes Gitter, rechts die ber\"uhmten 'Les Demoiselles d'Avignon' von Pablo Picasso (1907); sowohl die Relativit\"atstheorie als auch der Kubismus verbandeln auf ihre jeweilige Art Zeit und Raum miteinander. Mehr dazu bei Fischer \cite{epf}.}
%\end{figure}

Minkowskis {\it Raum} erstreckt sich mit seinen Wirkungsst\"atten K\"onigsberg, Bonn, Z\"urich und G\"ottingen \"uber weite Teile Europas und seine Leistungen in Mathematik und Physik waren bahnbrechend; seine {\it Zeit} hingegen war kurz: Hermann Minkowski verstirbt v\"ollig unerwartet an einem Blinddarmdurchbruch am 12. Januar 1909 in G\"ottingen. 

\section{Vorono\"\i{} und Blichfeldt}

Minkowskis Geometrie der Zahlen lieferte in relativ kurzer Zeit eine Vielzahl von neuen Resultaten. Es f\"allt auf, dass die nennenswerten hierunter w\"ahrend der Lebenszeit Minkowskis im Wesentlichen von ihm selbst stammen. Insofern mag man sich fragen, ob dies einzig auf seine Genialit\"at und Vertrautheit mit der neuen geometrischen Herangehensweise zur\"uckzuf\"uhren ist, oder seine Theorie zun\"achst nicht ordentlich rezipiert wurde.
\begin{quote}
''Die Puristen unter den Zahlentheoretikern, haben die Geometrie der Zahlen als der Analysis angeh\"orig gefunden und sich bem\"uht, dieses Werkzeug m\"oglichst zu eliminieren. So wurden f\"ur den Minkowski'schen Linearformensatz verschiedene, elementare Beweise gegeben. Am bekanntesten ist der Beweis von A. Hurwitz (1859-1919), der ihn zun\"achst f\"ur ganzzahlige Koeffizienten mit Hilfe des Schubfachprinzipes beweist. Daraus folgt unmittelbar der Fall der rationalen Koeffizienten und durch Grenz\"ubergang kann der allgemeine Fall erledigt werden.''
\end{quote}
schrieb Edmund Hlawka \cite{hlawka}, S. 406, selbst ein Forscher auf dem Gebiet der Geometrie der Zahlen. Explizit angesprochen wird hier der elementare Beweis \cite{hurwitzelem}, den Adolf Hurwitz, einer der pr\"agendsten Lehrer Minkowskis, aber auch dessen Freund und Kollege w\"ahrend seiner Z\"uricher Zeit, f\"ur dessen Linearformensatzes geliefert hatte. Nat\"urlich beleben verschiedene Ans\"atze und Beweise die Materie. Auch kann man den vielseitigen Hurwitz kaum als Puristen bezeichnen, der nicht aufgeschlossen f\"ur neue Methoden und Entwicklungen w\"are. Nicht zuletzt ist es aber die geometrische Intuition, welche die Minkowskische Geometrie der Zahlen so elegant erscheinen l\"asst und nebenbei derart tiefe arithmetische Resultate erlaubt. Tats\"achlich ist diese neue Denkweise nicht nur Minkowski eigen.      
\smallskip

Georgy Vorono\"\i{} hatte vieles mit Minkowski gemein (nicht nur sein Aussehen). Seine Lebensspanne deckt sich ziemlich mit der Minkowskis. Vorono\"\i{} wurde am 28. April 1868 in Zhuravka (damals Russland, heute Ukraine) geboren und verstarb jung und unerwartet am 20. November 1908 in Warschau an entz\"undeten Gallensteinen. Die beiden trafen sich ein einziges Mal und zwar auf dem Internationalen Mathematiker Kongress 1904 in Heidelberg, wo Vorono\"\i{} zum Kreisproblem vortrug \cite{voronoi04}, w\"ahrend Minkowski seine Geometrie der Zahlen vorstellte \cite{minkowski04}.
%\begin{figure}[ht]
%\centering
%\includegraphics[height=7.5cm]{minkowski2.eps}\quad \includegraphics[height=7.5cm]{voronoy.eps}
%\caption{\footnotesize Wer ist wer? Hermann Minkowski (links) und Georgy Vorono\"\i{} (rechts) sind nicht so leicht auseinanderzuhalten.}
%\end{figure}
Und Vorono\"\i{} gilt als Mitbegr\"under der Geometrie der Zahlen.\footnote{Je nach geographischem und kulturellem Hintergrund wird die Geometrie der Zahlen dem einen oder anderen zugeschrieben. Vielleicht ist die Geometrie der Zahlen ein weiteres Beispiel f\"ur das Ph\"anomen zeitgleicher unabh\"angiger Erkenntnisse wie etwa der Beweis des Primzahlsatzes.} 
\par

Bereits 1895 bestand ein Kontakt zwischen den beiden; so berichtet Minkowski in einem Brief vom 4. Dezember an seinen Freund Hilbert \cite{briefe}, S. 72, von einem 
\begin{quote}
''Woronoj (Petersb.) [der] auf 188 Seiten eine Theorie der kubischen K\"orper entwickelt, von der ich nach dem mangelhaften Referat in den Fortschritten jedenfalls nicht behaupten kann, da\ss{} sie die Hauptsachen au\ss er Acht lie\ss e. Du wirst in G\"ottingen vielleicht in der gl\"ucklichen Lage sein, einen das Russische verstehenden Mathematiker zur Verf\"ugung zu haben; dann vers\"aume es nicht, Dich \"uber jene Arbeiten zu informiren; sie k\"onnte doch mehr Gutes enthalten, als unsereiner von Russen erwartet.''
\end{quote}
Wir erinnern, dass Minkowskis Familie aufgrund von Antisemitismus aus dem russischen Aleksotas bzw. Kaunas ins preussische K\"onigsberg umsiedelte. Angesichts dessen zeugt dieses Zitat von einer verh\"altnism\"a\ss ig zur\"uckhaltenden Sicht, die keineswegs selbstverst\"andlich um die Jahrhundertwende war. Wie sich einem weiteren Brief vom 17. November 1896 entnehmen l\"asst, ergab sich im Weiteren sogar eine private Korrespondenz, die allerdings nicht erhalten geblieben ist. 
\par

Zu Beginn seiner Karriere besch\"aftigte sich Vorono\"\i{} mit Kettenbr\"uchen und deren Verwendung zur Konstruktion von Einheiten in Ganzheitsringen zu kubischen Z\"ahlk\"orpern; seine Doktorarbeit \cite{voronoiphd} erscheint 1896 nahezu zeitgleich mit Minkowskis verwandten \"Uberlegungen \cite{minko96} zu einer geometrischen Verallgemeinerung von Kettenbr\"uchen durch lokale Minima. Sp\"ater setzte Vorono\"\i{} sich zunehmend mit analytischen Themen zu divergenten Reihen und Mittelwerten zahlentheoretischer Funktionen auseinander. 
%\begin{figure}[ht]
%\centering
%\includegraphics[height=7.5cm]{vorono.eps}\quad \includegraphics[height=7.5cm]{postnihon.eps}%.eps}
%\caption{\footnotesize Links die Vorono\"\i{}-Zelle zu zuf\"allig gew\"ahlten Punkten der Ebene und rechts die Einzugsgebiete. Rechts eine Karte zum Auffinden der Postb\"uros in Tokio. Mit ein wenig Fantasie err\"at man, was hinter dem so genannten {\it post office problem} als Aufgabe steckt.}
%\end{figure}

Vorono\"\i{}s vielleicht wichtigster Beitrag zur Mathematik sind jedoch seine geometrischen \"Uberlegungen im Zusammenhang seiner Studien \cite{voronoi08} zu quadratischen Formen. In der euklidischen Ebene sei eine Menge ${\mathcal M}$ von Punkten gegeben; dann hei\ss t die Menge all der Punkte, welche von einem Punkt $P\in {\mathcal M}$ eine k\"urzere Distanz entfernt ist als von allen weiteren Punkten von ${\mathcal M}$, die {\it Vorono\"\i{}-Zelle} von $P$. Im Falle einer diskreten Menge ${\mathcal M}$ sind diese Zellen konvexe Polygone, welche \"uberschneidungsfrei den Raum pflastern. Dieses Konzept l\"asst sich leicht verallgemeinern (etwa auf andere R\"aume und andere Distanzma\ss e) und besitzt vielerlei Anwendungen, beispielsweise in der Meteorologie und Kristallographie\footnote{Tats\"achlich bemerkte bereits Seeber \cite{seeber} 1831, dass positiv definite quadratische Formen n\"utzlich f\"ur die Kristallographie sind.}, wo man auch gerne von Thiessen-Polygone sowie Wigner-Seitz-Zellen nach Alfred H. Thiessen (1911) bzw. Eugene Paul Wigner und Frederick Seitz (1933) spricht; tats\"achlich findet sich diese Idee bereits in einer Arbeit von Dirichlet \cite{dirichlet50} und noch fr\"uher in den {\it Le Monde, ou Traite de la lumiere} \cite{desca}, von Ren\'e Descartes in der Zeit um 1630 verfasst.
%\begin{figure}[ht]
%\centering
%\includegraphics[height=8.5cm]{descartes.eps}\quad \includegraphics[height=8.5cm]{descartespp.eps}
%\caption{\footnotesize Ren\'e Descartes (1596-1650) und eine Vorono\"\i{}-Zelle zwischen Sternen und Kometen}
%\end{figure}
Diese wichtige Serie von Arbeiten zu quadratischen Formen verfasste Vorono\"\i{} im Exil (\"ahnlich der Minkowskischen Isolation als mathematischer Pinguin in Bonn). Im Zuge der russischen Revolution 1905 und ihrer geographischen und zeitlichen Ausl\"aufer musste Vorono\"\i{} wie auch seine Kollegen bis ins Jahr 1908 (kurz vor seinem Tod) die geschlossene Universit\"at in Warschau, wo er seit 1894 t\"atig war, verlassen und lebte in dieser Zeit im fernen Novocherkassk in Russland. Dies erinnert an Newtons fruchtbaren R\"uckzug aus Cambridge angesichts der w\"utenden Pest 1665 und seine bahnbrechenden Gedanken in l\"andlicher Abgeschiedenheit. Mehr Informationen \"uber Vorono{\"{\i}}s kurzes Leben liefert Syta \cite{syta}. 
\par 

Einer der fr\"uhen Rezipienten und weiteren Ausgestalter der Minkowskischen Geometrie der Zahlen ist Hans Frederik Blichfeldt (1873-1945). Den geb\"urtigen D\"anen verschlug es als F\"unfzehnj\"ahrigen in die U.S.A., wo er zun\"achst mit harter Arbeit seinen Lebensunterhalt und das sp\"atere Studium finanzierte.\footnote{''I worked with my hands doing everything, East and West the country across.'' (cf. \cite{blichob}, S. 882)} Blichfeldts amerikanischer Traum wurde wahr --- er studierte erfolgreich in Stanford, erhielt dort eine Anstellung und promovierte dann 1898 bei Sophus Lie in Leipzig; sp\"ater wurde er Professor in Stanford. Seine Dissertation und auch seine fr\"uhen Arbeiten besch\"aftigen sich mit linearen Gruppen und deren Darstellungen, erst 1914 mit dem Artikel \cite{blichfeld} wendet er sich der Minkowskischen Ideenwelt zu und beweist den
\medskip

\noindent {\bf Satz von Blichfeldt (1914).} {\it Sei $\Lambda$ ein Gitter mit Determinate $\det(\Lambda)$ und $\Omega$ eine Punktemenge von einem Volumen ${\rm{vol}}(\Omega)> m\det(\Lambda)$ f\"ur eine nat\"urliche Zahl $m$. Dann existieren $m+1$ verschiedene Punkte ${\bf x}_0,\ldots,{\bf x}_m$ in $\Omega$, so dass die Differenzen ${\bf x}_j-{\bf x}_i$ allesamt in $\Lambda$ liegen.}
\medskip

\noindent Bemerkenswert ist hier die Losl\"osung von symmetrischen konvexen Mengen.\footnote{Ehrhart \cite{ehrhart55,ehrhart55b} entfernte f\"ur die Ebene die Symmetriebedingung durch Annahme des Schwerpunktes im Ursprung; im Allgemeinen ist die Frage jedoch offen.} Nennen wir eine Menge $\Omega\subset \R^n$ {\it packbar}, wenn $({\bf x}+\Omega)\cap ({\bf y}+\Omega)=\emptyset$ f\"ur verschiedene ${\bf x},{\bf y}\in\Z^n$ gilt. Dann besagt der Blichfeldtsche Satz somit: {\it Ist $\Omega$ Lebesgue-me\ss bar und packbar, dann gilt ${\rm{vol}}(\Omega)\leq 1$.} Dies liefert sofort einen alternativen Beweis des Minkowskischen Gitterpunktsatzes und sein Ansatz erlaubt dar\"uberhinaus noch weitere Anwendungen und Verallgemeinerungen. Einen interessanten alternativen Beweis des Blichfeldtschen Satzes mit Hilfe des Schubfachprinzips fand Willy Scherrer \cite{scherrer}.\footnote{Auch findet sich hier die Approximationsidee aus Mordells Beweis des Gitterpunktsatzes wieder.}
%\begin{figure}[ht]
%\centering
%\includegraphics[height=7cm]{blichfeldt.eps}\quad \includegraphics[height=7cm]{remak.eps}
%\caption{\footnotesize Von links nach rechts: Hans Frederik Blichfeldt (1873-1945), Robert Erich Remak (1888-1942)}
%\end{figure}

Eine weitere wesentliche Vereinfachung lieferte Robert Remak \cite{remak}. Der Berliner Remak promovierte 1911 bei Frobenius zu gruppentheoretischen Fragestellungen; seiner Habilitation wurde ihm jedoch unm\"oglich gemacht. Zwar hatte er bei seinem ersten Versuch 1919 u.a. den Frobenius' Nachfolger und Minkowski-Sch\"uler Constantin Carath\'eodory auf seiner Seite, aber seine schwierige Pers\"onlichkeit und offene Kritik am Lehrbetrieb standen seinen mathematischen Leistungen im Weg. Dar\"uberhinaus passten seine politischen Ansichten links der Mitte und auch seine Versuche, \"okonomische und soziale Theorien wissenschaftlicher zu fassen, nicht in das etablierte System der Berliner Universit\"at. Ein weiterer Versuch der Habilitation scheiterte 1923. Anschlie\ss end verbrachte Remak viel Zeit in G\"ottingen, gef\"ordert von Edmund Landau und Issai Schur, angeregt von Emmy Noether und -- erstaunlicherweise auch von dem Antisemiten Ludwig Bieberbach. In dieser Zeit forschte Remak zu Themen der Geometrie der 
Zahlen, fand Verallgemeinerungen der Minkowskischen S\"atze und etablierte sich damit schlie\ss lich derart, dass eine Habilitation in Berlin 1929 nicht mehr abgewiesen werden konnte. W\"ahrend der Nazizeit wurde Remak als Jude verfolgt und in Auschwitz ermordet.\footnote{Mehr Details zu Remaks Leben und Werk liefert Merzbach \cite{merzbach}. Beispielsweise lernt man dort, dass Remaks Gro\ss vater, ein bedeutender Mediziner gleichen Namens, als einer der Ersten seines Faches die au\ss erordentliche Erlaubnis erhielt, sich als Jude sich an der Berliner Universit\"at habilitieren zu d\"urfen, wurde aber trotz seiner Erfolge in Physiologie und Embryologie diskriminiert und blieb zeitlebens ohne eine ordentliche Professur.}
\par

\cite{minkowski01} hatte sich auch mit einer inhomogenen Version seines Linearformensatzes auseinandergesetzt und im Falle des Produktes von $n=2$ solcher Linearformen eine bestm\"ogliche Absch\"atzung einer minimalen nicht-trivialen L\"osung erzielt und eine Vermutung f\"ur allgemeines $n$ aufgestellt. Remak \cite{remak23} gelang mit erheblichem Aufwand unter Verwendung einer neuen Methode der Beweis einer solchen Absch\"atzung f\"ur den Fall $n=3$. Sp\"ater lieferten Freeman Dyson \cite{dys} (auf Anregung von Davenport, von dem sogleich die Rede sein wird) und Curtis McMullen \cite{mcm} Beweise f\"ur $n\leq 6$; weitere F\"alle wurden erfolgreich behandelt, worauf wir hier aber nicht weiter eingehen wollen.

\section{Die Schulen in Manchester und Wien}

Die nachfolgende zahlentheoretische Community der 1920er Jahre begegnete der Geometrie der Zahlen mit einer gewissen Zur\"uckhaltung, wie es ein weiteres Zitat von Hlawka belegt:\footnote{\cite{hlawka}, S. 404}
\begin{quote}
''Viele prominente Zahlentheoretiker (so z.B. H. Hasse (1898 geb.) waren der Ansicht, da\ss{} die Geometrie der Zahlen doch ein schwaches Werkzeug f\"ur die Zahlentheorie ist (...)''
\end{quote}
Erst eine weitere Generation von damals jungen und heute namhaften Vertretern der Zahlentheorie f\"uhrte zu einer angemessenen Rezeption. 
%\begin{figure}[ht]
%\centering
%\includegraphics[height=4.5cm]{mordell.eps}\quad \includegraphics[height=4.5cm]{mahler.eps}\quad \includegraphics[height=4.5cm]{davenport.eps}
%\caption{\footnotesize Wieder \"Ahnlichkeiten - von links nach rechts:  Louis Joel Mordell (1888-1972), dem man nachsagt, er habe an mehr als 170 verschiedenen Orten gelehrt und f\"ur dessen Biographie wir auf Cassels \cite{casselsmor} verweisen. In der Mitte: Kurt Mahler (1903-1988), Harold Davenport (1907-1969)}
%\end{figure}
Hlawka spricht von der Manchester Schule um Mordell und der Wiener Schule. 
\par

Behandeln wir zun\"achst die englische Schule. Sie wurde begr\"undet durch Louis Joel Mordell (1888-1972), Sohn j\"udischer Einwanderer aus Litauen (also eine gewisse geographische Verwandtschaft zu Minkowski); seit 1920 war er in Manchester t\"atig und baute dort eine bemerkenswerte Arbeitsgruppe auf.\footnote{Mordell schreibt ``I was self-taught mathematicially'' (cf. \cite{casselsmor}, S. 70) und tats\"achlich sagt man seinen Artikeln und B\"uchern einen eigent\"umlichen Stil nach.} Im Rahmen der Geometrie der Zahlen, welche Mordell in den 1930ern zunehmend besch\"aftigte, besteht sein Verdienst u.a. darin, erstmalig nicht-konvexe Probleme behandelt zu haben \cite{mordell41}. Der Austausch innerhalb der Manchester-Gruppe war sehr fruchtbar. Beispielsweise erzielte Mordells Sch\"uler Harold Davenport (wenngleich dieser bei John E. Littlewood promovierte) die exakten Minima f\"ur das Produkt dreier tern\"arer Linearformen \cite{hd38} und Mordell stellte dessen komplizierten Beweis einen eleganten Zugang zur Seite \cite{mordell42}. Kurze Zeit sp\"ater entdeckte Mordell \cite{mordell44} folgende sch\"one Abh\"angigkeit f\"ur aufeinanderfolgende Hermite-Konstanten:
$$
\gamma_n\leq \gamma_{n-1}^{(n-1)(n-2)}.
$$
Mordells Arbeitsgruppe wurde bereits in den fr\"uhen 1930ern Anlaufstelle f\"ur etliche aus Nazideutschland emigrierte Mathematiker, wie etwa Kurt Mahler, ein Sch\"uler Carl Ludwig Siegels j\"udischer Herkunft. Mahler initiierte u.a. die Behandlung von Sterngebieten\footnote{die sich dadurch auszeichnen, dass sie einen ausgezeichneten Punkt enthalten, so dass jede geradlinige Verbindungsstrecke zu einem weiteren Punkt komplett innerhalb dieser Menge liegt} mit Methoden der Geometrie der Zahlen. Sein Kompaktheitssatz \cite{mahler45} ist von konzeptioneller Bedeutung und ein wichtiges Werkzeug in Gregori Aleksandrovich Margulis' Beweis \cite{margulis87}, dass jede nicht-degenerierte quadratische Form $Q$ in $n\geq 3$ Ver\"anderlichen, welche nicht ein Vielfaches einer rationalen Form ist, ein dichtes Bild $Q(\Z^n)$ in $\R$ besitzt. Dies verifiziert eine alte Vermutung von Alexander Oppenheim \cite{oppenheim} bzw. einer versch\"arften Formulierung von Davenport.

%\begin{figure}[ht]
%\centering
%\includegraphics[height=5cm]{rogers.eps}\quad \includegraphics[height=5cm]{corput.eps}
%\caption{\footnotesize Links: Claude Ambrose Rogers (1920-2005); rechts: Johannes Gualtherus van der Corput (1890-1975). Letztgenannter spielt mit seinen Arbeiten zu Exponentialsummen eine gro\ss e Rolle in der Theorie der Zetafunktion und der Gleichverteilungslehre.}
%\end{figure}

Es sind weitere Charaktere an dieser Stelle zu nennen. Zum einen Claude Ambrose Rogers, welcher bereits w\"ahrend seiner Promotion mit etlichen Ver\"offentlichungen zur Geometrie der Zahlen auffiel und in dieser Zeit viel mit Davenport kooperierte. Zum anderen Johannes van der Corput \cite{vdc}, der nicht direkt der Manchester Schule zugeordnet werden kann, sich daf\"ur aber als Doktorvater von Jurjen Koksma auszeichnete. Gemeinsam mit Davenport \cite{vdchd} gelang eine Verallgemeinerung des Minkowskischen Gitterpunktsatzes in der Ebene, indem die Kr\"ummung der Randkurve Ber\"ucksichtigung findet. Ferner entwickelte van der Corput eine Methode Mordells \cite{mordell34} (welche auch den Beweis des Minkowskischen Gitterpunktsatzes aus \S 3 liefert) fort, welche u.a. eine untere Absch\"atzung f\"ur die Anzahl der Gitterpunkte in einem symmetrischen konvexen K\"orper in Abh\"angigkeit von dessen Volumen erlaubt. Eine verwandte Verallgemeinerung von Siegel \cite{siegel35} mittels Fourier-Analysis brachte nicht 
nur neue analytische Werkzeuge ins Spiel, sondern lieferte sogar eine explizite (gewichtete) Gleichung.    
\medskip 

Die Wiener Schule wurde von Philipp Furtw\"angler\footnote{der an demselben Gymnasium Andreanum in Hidlesheim wie Adolf Hurwitz zur Schule ging} gegr\"undet und sp\"ater von Nikolaus Hofreiter und nicht zuletzt Edmund Hlawka fortgef\"uhrt. 
%\begin{figure}[ht]
%\centering
%\includegraphics[height=5cm]{furtwaengler.eps}\quad \includegraphics[height=5cm]{hofreiter.eps}\quad \includegraphics[height=5cm]{hlawka.eps}
%\caption{\footnotesize Von links nach rechts: Friedrich Pius Philipp Furtw\"angler (1869-1940), Nikolaus Hofreiter (1904-1990) und Edmund Hlawka (1916-2009)}
%\end{figure}
Furtw\"angler ver\"offentlichte recht wenig zur Geometrie der Zahlen, gab allerdings viele Anst\"o\ss e, insbesondere f\"ur die Forschung seines Sch\"ulers Hofreiter, welcher sich intensiv mit Produkten von inhomogenen Linearformen, quadratischen Zahlk\"orpern ohne euklidischen Algorithmus und Fragen der diophantischen Approximation auseinandersetzte \cite{hofr}. In der Regel sind die Ganzheitsringe von Zahlk\"orpern nicht-euklidisch, womit zun\"achst kein euklidischer Algorithmus zur Verf\"ugung steht mit Hilfe dessen ein praktikabler Kettenbruchalgorithmus implementiert werden k\"onnte. Hier liefert der Minkowskische Linearformensatz einen alternativen Zugang.\footnote{Ein typisches Beispiel findet sich in der Dissertation von Hilde Gintner \cite{gintner} aus dem Jahr 1936: {\it Zu  $m\in\N$ und beliebigem $\alpha\in\C\setminus \Q(i\sqrt{m})$ existieren unendlich viele $p,q\in\Z[i\sqrt{m}]$ mit} $\vert \alpha-{p\over q}\vert<{\sqrt{6m}\over \pi}\,{1\over \vert q\vert^2}$. Im Falle $m\not\equiv 3\bmod\,4$ 
ist $\Z[i\sqrt{m}]$ der Ganzheitsring des Zahlk\"orpers $\Q(i\sqrt{m})$; andernfalls liefert eine Beweisvariante eine analoge Ungleichung f\"ur den entsprechenden Ganzheitsring als rechte Seite sogar ${\sqrt{6m}\over 2\pi}\,{1\over \vert q\vert^2}$.} Die relevanten Beitr\"age der Wiener Schule zur Geometrie der Zahlen sind aber Edmund Hlawka zuzuschreiben. In seiner Promotion \"uber die Approximationen von zwei komplexen inhomogenen Linearformen \cite{drhlawka} bei Hofreiter f\"uhrte er dessen Steckenpferd fort; im Rahmen seiner Habilitation gelang ihm dann ein gro\ss er Wurf, f\"ur dessen Erkl\"arung wir allerdings etwas ausholen m\"ussen. 
\par

F\"ur eine kompakte Menge ${\mathcal C}\subset \R^n$ ist die zugeh\"orige {\it Gitterkonstante} definiert als 
$$
\Delta({\mathcal C})=\min\{\det(\Lambda)\,:\, \Lambda\cap{\mathcal C}\neq\{{\bf 0}\}\},
$$
wobei das Minimum \"uber alle Gitter des $\R^n$ erhoben wird (und stets existiert wie Mahler zeigte). Dann ist 
$$
2^{-n}{{\rm vol}({\mathcal C})\over \Delta({\mathcal C})}
$$
gleich dem Maximum der {\it Gitterpackungsdichte}, also dem Anteil des $\R^n$, welcher durch ${\mathcal C}+\Lambda$ ohne \"Uberlappung \"uberdeckt wird. Im speziellen Fall von Kugeln notieren wir diese Dichte als $\delta_n$. Den Zusammenhang zur Hermite-Konstanten $\gamma_n$ hatte bereits Gau\ss{} \cite{gauss40} mit der Formel
$$
\delta_n={1\over \Gamma(1+{\textstyle{n\over 2}})}\left({\pi\gamma_n\over 4}\right)^{n\over 2}
$$
festgestellt. Eine untere Absch\"atzung f\"ur diese wichtige Gr\"o\ss e und damit auch f\"ur die Hermite-Konstante liefert der
\bigskip

\noindent {\bf Satz von Hlawka-Minkowski (1911, 1943).} {\it F\"ur einen Sternk\"orper ${\mathcal S}\subset \R^n$ gilt} 
$$
\Delta({\mathcal S})\leq {{\rm{vol}}({\mathcal S})\over 2\zeta(n)}.
$$ 
\medskip

\noindent Hierbei ist $\zeta(n)=1+2^{-n}+3^{-n}+\ldots$ der Wert der Riemannschen Zetafunktion $\zeta(s)$ an der Stelle $s=n$. Dieses Resultat findet sich ohne Beweis in den gesammelten Abhandlungen von Minkowski (posthum publiziert 1911); den ersten Beweis lieferte Hlawka \cite{hlawkahabil} in seiner Habilitationsschrift 1943. Kurze Zeit sp\"ater ver\"offentliche Siegel \cite{siegel45} einen weiteren, vielleicht noch zug\"anglicheren Beweis. Hlawka betreute Zeit seines langen Lebens 130 Promotionen und etliche seiner Sch\"uler haben Professuren inne, etliche von ihnen mit einem Forschungsgebiet nahe der Geometrie der Zahlen. ''Die st\"urmische Entwicklung dieses Gebietes in den 40- und 50-ziger Jahren endet ungef\"ahr um 1960. Das Lehrbuch von Lekkerkerker: 'Geometry of Numbers' (1969) stellt das Erzielte zusammen. Die Entwicklung ging aber in stilleren Bahnen weiter.'' schreibt Hlawka \cite{hlawka}, S. 9. Die bahnbrechende Arbeit \cite{schmi} des Hlawka-Sch\"ulers Wolfgang Schmidt f\"allt in die 
Schlussphase der st\"urmischen Entwicklung. Der Schmidtsche {\it Teilraum-Satz}\footnote{in der mittlerweile englischsprachigen Literatur `subspace theorem'} besagt, dass {\it zu gegebenen linear unabh\"angigen Linearformen $Y_1,\ldots,Y_n$ in $n$ Unbekannten $X_1,\ldots,X_n$ mit algebraischen Koeffizienten sowie einem beliebigen $\epsilon>0$ ganze Zahlen $x_1,\ldots,x_n$ gibt, nicht alle null, so dass
$$  
\vert Y_1({\bf x})\cdot\ldots\cdot Y_n({\bf x})\vert < \vert {\bf x}\vert^{-\epsilon}\qquad\mbox{mit}\quad {\bf x}=(x_1,\ldots,x_n)
$$
in endlich vielen echten Unterr\"aumen des $\Q^n$ liegen.} 

%\begin{figure}[ht]
%\centering
%\includegraphics[height=5cm]{Kepler.eps}\quad \includegraphics[height=5cm]{thue.eps}\quad \includegraphics[height=5cm]{fejestoth.eps}
%\caption{\footnotesize Links der Urheber eines gro\ss en Problems: Johannes Kepler (1571-1630); in der Mitte ein Probleml\"oser: Axel Thue (1863-1922); rechts ein Wegbereiter: L\'aszl\'o Fejes T\'oth (1915-2005)}
%\end{figure}

\section{Packungsprobleme}

Nicht unerw\"ahnt bleiben d\"urfen die so genannten Packungsprobleme. Hierbei geht es weniger um Fragen der Optimierung wie sie Paketzustelldienste oder moderne B\"ucherverkaufskonzerne zu l\"osen haben, als um eine naheliegende Frage, die sich bereits zu Beginn dieses Artikels angesichts des Zusammenhangs von Dreieckszahlen und Kugelpyramiden oder aber beim Besuch eines Wochenmarktes stellen mag: Wie lassen sich m\"oglichst viele Orangen auf engem Raum stapeln? 
\par

Diese Fragestellung ist tats\"achlich recht alt. Der englische Astronom Thomas Harriot war nicht nur 1609 der Erste (vor Galilei!), der sein Fernrohr gen Himmel neigte und dort u.a. die Jupitermonde entdeckte, was er jedoch mitzuteilen unterlie\ss{}, sondern er machte sich auch Gedanken \"uber das raumsparendste Stapeln von Kanonenkugeln (auf Anregung des auf den Weltmeeren herumsegelnden Sir Walter Rayleigh).\footnote{Wesentlich \"alter sind Untersuchungen von ArybhatIya und Bhaskhara aus dem f\"unften bzw. sechsten Jahrhundert in Sanskrittexten; cf. Hales \cite{haleshist}, S. 2.} Harriot \"au\ss erte die Vermutung, dass ein Stapeln derselben, wie man es von etwa Orangen auf einem Wochenmarkt kennt, optimal seien. Er war hierbei sicherlich nicht der Erste, denn die Orangenstapel entstehen ja nicht von ungef\"ahr. Ferner gab seine Korrespondenz mit dem noch bekannteren Astronomen und Mathematiker Johannes Kepler Anlass zu dessen \"Au\ss erung, dass {\it die dichteste Kugelpackung im dreidimensionalen 
euklidischen Raum durch eine kubisch-fl\"achenzentrierte Packung bzw. die hexagonale Packung gegeben ist.} Bemerkenswert ist hier die Losl\"osung von einem endlichen zu einem unendlichen Orangenstapel.
%\begin{figure}[ht]
%\centering
%\includegraphics[height=7cm]{keplerconjecture.eps}
%\caption{\footnotesize Wo ist der Unterschied? Es besteht eine gro\ss e \"Ahnlichkeit zwischen der kubisch-fl\"achenzentrierten Kugelgitterpackung und der hexagonalen Kugelgitterpackung; die zugeh\"orige Packungsdichte betr\"agt in beiden F\"allen optimale $74,048\ldots$ Prozent.}
%\end{figure}
Diese so genannte Keplersche Vermutung findet sich in dessen Arbeit {\it vom sechseckigen Schnee} \cite{kepler} von 1610. Aus den romantischen Eingangsworten Keplers mit der Widmung f\"ur seinen Freund Johannes Wackher von Wackersfeld mag man den Zusammenhang mit der Kugelpackung heraus erahnen:\footnote{siehe auch {\sf http://www.mathematik.de/ger/information/kalenderblatt/keplers\_strena/keplers\_strena.html}} 
\begin{quote}
''Als ich so nachdenklich und sorgenvoll \"uber die Br\"ucke ging und mich \"uber meine eigene Armseligkeit \"argerte, n\"amlich zu dir ohne Neujahrsgeschenk zu kommen, und immer denselben Gedanken nachging, dieses Nichts anzugeben oder etwas zu finden, was ihm am n\"achsten kommt, und ich daran die Sch\"arfe meines Denkens \"ubte, da f\"ugte es der Zufall, dass sich der Wasserdampf durch die K\"alte zu Schnee verdichtete und vereinzelte kleine Flocken auf meinen Rock fielen, alle waren sechseckig mit gefiederten Strahlen. [\ldots] Ei, das ist ein erw\"unschtes Neujahrsgeschenk f\"ur einen Freund des Nichts! So wie der Schnee da vom Himmel herabkommt und den Sternen \"ahnlich ist, ist er auch passend als Geschenk eines Mathematikers, der nichts hat und nichts erh\"alt.'' 
\end{quote}
Obwohl zu Keplers Zeiten die atomare bzw. molekulare Struktur unbekannt ist, bemerkt dieser eine hexagonale Anordnung (welche die Wassermolek\"ule bilden). Diesen Gedanken wird Auguste Bravais 1849 mit der Idee seiner Kristallgitter theoretisch genauer fassen und Max von der Laue wird diese 1912 experimentell nachweisen.\footnote{Der Nachweis der Atomstruktur schloss sich an die 1905 von Einstein untersuchten Brownschen Bewegung an, welche auch in der Mathematik eine wichtige Rolle in der Theorie der stochastischen Prozesse spielt. Die Idee der {\it unteilbaren} Materie findet sich bereits bei Demokrit, und die Primzahlen sind deren mathematische Entsprechung.} Diese augenscheinliche Struktur der Schneeflocken erinnerte Kepler an seine Konversation mit Harriot und f\"uhrte ihn zu seiner Vermutung zur raumsparenden Lagerung unendlich vieler Kugeln gleicher Gr\"o\ss e.  
%\begin{figure}[ht]
%\centering
%\includegraphics[height=6cm]{kissingnumber2.eps}\qquad \includegraphics[height=6cm]{Kissing.eps}
%\caption{\footnotesize Die Kusszahlen f\"ur die Dimensionen $n=2$ und $3$ sind $6$ bzw. $12$. Bereits f\"ur den $10$-dimensionalen Raum ist die Kusszahl unbekannt.}
%\end{figure}

Die analoge Frage f\"ur Kreise ist deutlich einfacher, auch liefert sie ein wenig \"Uberzeugung f\"ur die Richtigkeit der Keplerschen Vermutung. Ein wenig ebene Geometrie bzw. Herumexperimentieren mit einem Zirkel legt nahe, dass man um einen Kreis sechs weitere gleich gro\ss e Kreise finden kann, die diesen ber\"uhren und sich dabei nicht \"uberlappen. Entsprechend ist $6$ die {\it Kusszahl}\footnote{im Englischen `kissing number'} in der Ebene. F\"ur den dreidimensionalen Raum ist die Situation verzwickter: Hier stritten Isaac Newton und David Gregory Ende des 17. Jahrunderts, ob es $12$ oder gar $13$ Kugeln sind, welche eine gegebene Kugel gleicher Gr\"o\ss e ohne \"Uberlappung ber\"uhren k\"onnen. Tats\"achlich hatte Newton mit seiner $12$ recht, wenngleich Platz f\"ur eine dreizehnte Kugel best\"unde, nicht jedoch als Ganzes, sondern lediglich zerchnitten. Der erste rigorose Beweis hierf\"ur stammt von Kurt Sch\"utte \& Bartel Leendert van der Waerden \cite{53} aus dem zwanzigsten Jahrhundert.\footnote{
L\'aszl\'o Fejes T\'oth \cite{toth} hat weitere Fragen zu endlichen Packungen aufgeworfen; diese stellen insbesondere mit Blick auf au\ss ermathematische Anwendungen eine Herausforderung dar. Leppmeier \cite{leppmeier} behandelt u.a diese f\"ur eine Vielzahl von Optimierungsproblemen interessanten Fragestellungen (wie etwa platzsparende Packungen von Tennisb\"allen) und auftretende \"Uberraschungen, wof\"ur an dieser Stelle lediglich das Stichwort {\it Wurstkatastrophe} genannt sei.}     
\par

F\"ur das Folgende ist es sinnvoll, einige Definitionen zu treffen. Unter einer {\it Kugelpackung} $\kappa$ sei ein Arrangement von unendlich vielen gleichgro\ss en Kugeln verstanden, die sich nur ber\"uhren, aber nicht in mehr als einem Punkt \"uberlappen d\"urfen. Die {\it Dichte} $\delta(\kappa)$ einer solchen Kugelpackung definieren wir als
$$
\delta(\kappa)=\lim_{r\to\infty}{{\rm{vol}}(\kappa\cap {\mathcal B}_r)\over {\rm{vol}}({\mathcal B}_r)},
$$
wobei ${\mathcal B}_r$ eine Kugel vom Radius $r$ im $\R^n$ bezeichne. Damit ist $\delta(\kappa)$ der {\it Anteil} von $\kappa$, der in ${\mathcal B}_r$ im Grenzwert $r\to\infty$ enthalten ist. Wir sprechen von einer {\it Kugelgitterpackung}, wenn die Mittelpunkte der Kugeln ein Gitter bilden; ferner sprechen wir in der Ebene von Kreisen anstelle von Kugeln.
\par

Die dichteste Kreisgitterpackung ist, wie Lagrange \cite{lagrange} 1773 zeigte, hexagonal und bedeckt ${\pi\over 2\sqrt{3}}=90,69\ldots$ Prozent der Ebene; es gilt also
$$
{\pi\over 2\sqrt{3}}=\max_\kappa \delta(\kappa)
$$
und die Dichte ist maximal unter aller Kreisgitterpackungen $\kappa$ f\"ur ein hexagonales Gitter; in diesem Zusammenhang liefert die bin\"are quadratische Form $X^2+XY+Y^2$ das zugeh\"orige Minimum\footnote{und die Eisensteinschen Zahlen $m+{1\over 2}(-1+\sqrt{-3})n$ mit ganzzahligen $m,n$ bilden dieses Gitter in der komplexen Ebene mittels der dritten Einheitswurzel ${1\over 2}(-1+\sqrt{-3})=\exp({2\pi i\over 3})$.}. Diese Extremaleigenschaft suggeriert bereits die Kusszahl im Zweidimensionalen, aber dies ist nat\"urlich kein Beweis. Die dichteste Kugelgitterpackung im dreidimensionalen Raum ist, wie Gau\ss{} \cite{gauss40} 1831 herleitete, fl\"achenzentriert-kubisch mit einem Bedeckungsanteil von ${\pi\over 3\sqrt{2}}=74,04\ldots$ Prozent; die zugeh\"orige tern\"are quadratische Form ist $X^2+Y^2+Z^2+XY+YZ+XZ$. Tats\"achlich sind beide Packungen sogar unter allen Kreis- bzw. Kugelpackungen {\it optimal}, wie Orangenverk\"aufer und sogar mathematisch ungebildete Bienenv\"olker wissen. Aber hierf\"ur einen rigorosen mathematischen Beweis zu geben, erwies sich klange Jahre als extrem schwierig. Die Abwesenheit von Struktur ist hierf\"ur der Grund: Liegt eine Gitterstruktur vor, so l\"asst sich die damit verbundene Regelm\"a\ss igkeit in der Anordnung relativ einfach nutzen; hingegen ist eine aperiodische Packung wesentlich unangenehmer zu untersuchen. Das zweidimensionale Analogon der Keplerschen Vermutung behandelte Axel Thue \cite{thue92,thue10} erfolgreich. Seine Idee benutzt implizit das Konzept der Vorono{\"{\i}}-Zellen. Diese erweisen sich im Falle der dichtesten Kreispackung der Ebene als regul\"are Sechsecke, womit folglich das hexagonale Gitte auf die dichteste Kreispackung f\"uhrt -- ein Grund mehr, auch f\"ur den dreidimensionalen Fall zu erwarten, dass die dichteste Kugelpackung von einem Gitter herr\"uhrt.    
%\begin{figure}[ht]
%\includegraphics[height=6.5cm]{bienenwaben.eps}
%\caption{\footnotesize Bienen beim Bau ihrer stabilen und dabei optimierten Waben mittels einer hexagonalen Anordnung}
%\end{figure}

Die Keplersche Vermutung wurde von Hilbert in seine Liste der mathematischen Probleme f\"ur das zwanzigste Jahrhundert (als Teil des allgemeiner formulierten 18. Problems) aufgenommen (siehe \cite{hilberticm}). Aber erst vor kurzem gelang es Thomas Callister Hales diese Nuss mit massiven Computereinsatz zu knacken.\footnote{Sein \"alterer Namensvetter, der Physiologe Stephen Hales (1677-1761) stellte im Rahmen seiner Untersuchungen zur Pflanzenphysiologie \cite{shales} randomisierte Experimente mit Erbsen an, welche f\"ur die Keplersche Vermutung sprachen.} Nach mehrj\"ahriger Arbeit reichte er 1998 seine Arbeit in den renommierten {\it Annals of Mathematics} ein. Sein Ansatz basierte auf der oben angesprochenen alten Idee von L\'aszl\'o Fejes T\'oth\footnote{Hales listet aber noch weitere Unterst\"utzer auf; insbesondere die Arbeit von Samuel P. Ferguson ist hier zu nennen.} \cite{toth} aus dem Jahre 1953, welche an Thues Ansatz ankn\"upfte und es erlaubt, die potentiell unendlich vielen Kugelanordnungen 
auf eine endliche Anzahl von F\"allen zu reduzieren; wesentliches Werkzeug dabei sind die Vorono{\"{\i}}-Zellen. Schlie\ss lich war Hales nach jahrelanger Vorarbeit soweit, dass lediglich einhundert F\"alle \"ubrigbleiben, die er dann einzeln mit den schnellsten verf\"ugbaren Rechnern behandelt.\footnote{Dies erinnert an den Computerbeweis des Vierfarbensatzes durch Kenneth Apel \& Wolfgang Haken 1976.} Aber gerade dieser Einsatz von Computern gibt letztlich Anlass zur Skepsis. Hales' umfassender Beweis von ca. 250 Seiten und 40\,000 Zeilen Computercode wurde von Dutzenden Gutachtern \"uber Jahre hinweg analysiert. Schlie\ss lich wird Hales Arbeit \cite{hales0} 2005 ver\"offentlicht, allerdings mit dem Zusatz, dass die Gutachter zu 99 Prozent von der Richtigkeit \"uberzeugt seien. Hier schwingt auch Unmut mit: Ein Beweis soll einen mathematischen Sachverhalt erleuchten; gewisserma\ss en ist der Beweis wichtiger als die Aussage selbst! Und dies geht Hales' Computerbewies ab. Nun haben sich Mathematik und 
Informatik seitdem rasant weiter entwickelt. Mittlerweile existieren so genannte {\it Beweisassistenten}, die im Gegensatz zu herk\"ommlichen Computerprogrammen komplexe Rechnungen gem\"a\ss{} den Reglen der Logik verifizieren. Hales m\"ochte den Makel, der seinem Beweis anhaftete, mit einem solchen Beweisassistenten beseitigen. Er initiiert ein Projekt mit Namen {\it Formal Proof of Kepler}, kurz {\it flyspeck}\footnote{also `Fliegendreck'}, welches mit staatlicher Unterst\"utzung und Einsatz des Informatikers Georges Gonthier von Microsoft tats\"achlich in erstaunlicher Geschwindigkeit innerhalb eines Monats 2014 ein zweites Mal die Keplersche Vermutung verifiziert.
%\begin{figure}[ht]
%\centering
%\includegraphics[height=5cm]{thcahales.eps}\quad \includegraphics[height=5cm]{keplerverm.eps}
%\caption{\footnotesize Links: Thomas Callister Hales (1958-\ ) und sein Orangenstapel; das Bild entstammt der Zeit-Ausgabe, 2. Juli 2015 mit dem lesenswerten Artikel 'Das Orangen-Projekt' von Rudy Novotny. Auch das Buch \cite{szpiro} von George Szpiro zur Keplerschen Vermutung sei an dieser Stelle empfohlen; hier lernt die Leser\_in nicht nur einiges zu Packungsproblemen und deren bewegter Geschichte, aber auch anderes Wissenswertes, wie etwa die Antwort auf die Frage, warum die Natur so viele nahezu vollkommen runde Fr\"uchte und Gem\"use hervorbringt. Rechts: Eine Illustration aus Keplers `Strena seu de Nive sexangula', die seine nun bewiesene Vermutung illustriert.}
%\end{figure}

In fast allen h\"oheren Dimensionen ist die Frage nach der dichtesten Kugelpackung, ja sogar nach der dichtesten Kugelgitterpackung bislang ungel\"ost. Lediglich f\"ur Spezialf\"alle sind Antworten bekannt. Eine besondere Rolle spielt hier das Leech-Gitter, entdeckt 1967 von John Leech \cite{leech}, welches im $24$-dimensionalen Raum die dichteste Kugelpackung bereitstellt und auf die beeindruckende Kusszahl $196\,560$ f\"uhrt, wie Henry Cohn et al. \cite{via2} k\"urzlich zeigten; zuvor hatte Maryna Viazovska \cite{via1} den verwandten Fall der Dimension $8$ gel\"ost mittels eines Zusammenhangs mit Quasimodulformen. Die speziellen Eigenschaften des Leech-Gitters machen es in vielerlei Weise interessant f\"ur Algebra und Zahlentheorie (z.B. hinsichtlich der sporadischen Gruppen und Modulformen), aber auch in der Codierungstheorie finden sich Anwendungen, worauf wir hier aber nicht weiter eingehen wollen (sondern stattdessen auf das Standardwerk von Conway \& Sloane \cite{conway} verweisen). Erstaunlicherweise wird in allgemeinen h\"oheren Dimensionen erwartet, dass die dichtesten Kugelpackungen nicht von Gittern herr\"uhren, da die Struktur von Gittern zu einschr\"ankend ist (cf. Cohn \& Elkies \cite{cohnelkies}, S. 690).

Minkowski selbst sah seine Geometrie der Zahlen als Wegbereiter:\footnote{\cite{minkowski07}, S. 234/235}
\begin{quote}
''Alle Theoreme hier wiesen {\it einen} Ursprung auf, wir sch\"opften sie aus einer gemeinsamen, sehr durchsichtigen Quelle, die ich als das {\it Prinzip der zentrierten konvexen K\"orper im Zahlengitter} bezeichnen m\"ochte. Nun sind wir in der Tat eine Strecke Wegs in das Reich der heutigen Zahlentheorie eingedrungen. Wir k\"onnen daran denken, uns auf diesem Boden zu akklimatisieren.'' 
\end{quote}
Anschlie\ss end formulierte Minkowski, dass das Studium der Primideale wunderbare Zusammenh\"ange zwischen der Zahlentheorie und der Theorie der Funktionen offenbaren w\"urde. Tats\"achlich sind auf den Gebieten der algebraischen Zahlentheorie und in der Primzahlverteilung in den fr\"uhen Jahren des zwanzigsten Jahrhunderts erhebliche Fortschritte erzielt worden.

\medskip

\hspace*{3.7cm} \hrulefill 

\bigskip 

\noindent Einige Aspekte der Geometrie der Zahlen konnten nicht in ihrer vollst\"andigen Tiefe untersucht werden. Beispielsweise sei die additive Kombinatorik von Imre Ruzsa \cite{ruzsa} erw\"ahnt, in der Fragen nach dem Pendant von Ma\ss en und Dimension bei diskreten Mengen nachgegangen wird. F\"ur die von Eug\`ene Ehrhart \cite{ehrhart62} angestossenen umfangreichen Untersuchungen und Verallgemeinerungen des Pickschen Satzes, welche Geometrie, Kombinatorik und Zahlentheorie auf sehr sch\"one Weise miteinander verkn\"upfen, sei hier auf das sch\"one Buch von Beck \& Robins \cite{beck} verwiesen; Hugo Hadwigers Arbeiten zur Zerlegungen von Polytopen, stimuliert durch Max Dehns L¿``osung des dritten Hilbertschen Problems, besprechen Erd\H{o}s et al. \cite{erdoes}. Verallgemeinerungen der Minkowskischen Theorie f\"ur Zahlk\"orper werden von K. Rogers \& P. Swinnerton-Dyer \cite{krogers} behandelt. Interessant sind auch Anwendungen der Geometrie der Zahlen und der Theorie der Gitter auf Fragestellungen der Kodierungstheorie und der Kryptographie (wie etwa das Auffinden eines k\"urzesten Gittervektors). Hier lieferten Lenstra, Lenstra \& Lov\'asz \cite{lll} mit ihrem LLL-Algorithmus zur Bestimmung kurzer Gittervektoren in Polynomialzeit einen \"au\ss erst wichtigen Beitrag (und eine Verbesserung und Erweiterung eines klassischen Verfahrens von Gau\ss ). F\"ur Anwendungen der Geometrie der Zahlen au\ss erhalb der Zahlentheorie und sogar der Mathematik verweisen wir auf Lov\'asz \cite{lovasz}. F\"ur ein tieferes Eindringen in die Geometrie der Zahlen empfehlen wir die B\"ucher von Gruber \& Lekkerkerker \cite{gruber} sowie die Klassiker von Cassels \cite{cassels} und Siegel \cite{siegel}; geschichtliche Aspekte insbesondere in einem gr\"o\ss eren geometrischen Kontext findet man bei Gruber \cite{gruber90} und Opolka \& Scharlau \cite{scharl}. F\"ur die Gitterpunktprobleme bietet das Buch von Fricker \cite{fricker} einen sehr guten Einstieg. Packungsprobleme werden sehr sch\"on bei Rogers \cite{rogers} dargestellt. Eine detaillierte Geschichte der Geometrie der Zahlen liefert die Dissertation von S\'ebastien Gauthier \cite{gauthier}.

\bigskip

\noindent Nicola Oswald, Institut f\"ur Mathematik und Informatik, Bergische Universit\"at Wuppertal, Gau\ss str. 20, 42\,119 Wuppertal, Germany, oswald@uni-wuppertal.de;\\
%\quad und\\
%Institut f\"ur Mathematik, Julius-Maximilians Universit\"at W\"urzburg, Emil-Fischer-Str. 40, 97\,074 W\"urzburg, Germany, nicola.oswald@mathematik.uni-wuerzburg.de

%\newpage

%\begin{figure}[ht]
%\includegraphics[height=13cm,angle=90]{europa1896.eps}
%\caption{\footnotesize Neuer Handatlas, E. Debes (Herausgeb.), Wagner \& Debes, Leipzig 1896, Karte Nr. 11, Europa. Mit Blick auf die Protagonisten sei auf dieser zeitgen\"ossischen Karte auf deren Lebensmittelpunkte in Bonn, Z\"urich, G\"ottingen und Warschau hingewiesen.}
%\end{figure}

\end{document}